\documentclass[11pt]{amsart}
\reversemarginpar
\usepackage{amsmath} 
\usepackage{amssymb}
\usepackage{euscript}

\newcommand{\nwc}{\newcommand}

\nwc{\nwt}{\newtheorem}

\nwt{prop}{Proposition}
\nwt{proposition}{Proposition}
\nwt{lem}{Lemma}
\nwt{thm}{Theorem}
\nwt{cor}{Corollary}
\nwt{rem}{Remark}
\nwt{defin}{Definition} 

\nwc{\barr}{\begin{array}}
\nwc{\bal}{\begin{align}}
\nwc{\bequ}{\begin{equation}}
\nwc{\ben}{\begin{equation*}}
\nwc{\bea}{\begin{eqnarray}}
\nwc{\beq}{\begin{eqnarray}}
\nwc{\bean}{\begin{eqnarray*}}
\nwc{\beqn}{\begin{eqnarray*}}
\nwc{\beqast}{\begin{eqnarray*}}

\nwc{\earr}{\end{array}}
\nwc{\eal}{\end{align}}
\nwc{\eequ}{\end{equation}}
\nwc{\een}{\end{equation*}}
\nwc{\eea}{\end{eqnarray}}
\nwc{\eeq}{\end{eqnarray}}
\nwc{\eean}{\end{eqnarray*}}
\nwc{\eeqn}{\end{eqnarray*}}
\nwc{\eeqast}{\end{eqnarray*}}

\nwc{\bldm}{\boldmath} 
\nwc{\ubm}{\unboldmath} 
\nwc{\mf}{\mathbf} 
\nwc{\blds}{\boldsymbol} 
\nwc{\ml}{\mathcal} 

\nwc{\al}{\alpha}
\nwc{\vep}{\varepsilon}
\nwc{\ep}{\epsilon}
\nwc{\veps}{\varepsilon}
\nwc{\eps}{\epsilon}
\nwc{\vrho}{\varrho}
\nwc{\orho}{\bar\varrho}
\nwc{\vpsi}{\varpsi}
\nwc{\lamb}{\lambda_\varepsilon}
\nwc{\lam}{\lambda}
\nwc{\del}{\delta}
\nwc{\tht}{\theta}
\nwc{\om}{\omega}

\nwc{\btht}{\blds{\theta}}
\nwc{\bxi}{\blds{\xi}}
\nwc{\bmal}{\blds{\al}}
\nwc{\bmbe}{\blds{\beta}}

\nwc{\la}{\langle}    
\nwc{\ra}{\rangle}    
\nwc{\lp}{\left(}     
\nwc{\rp}{\right)}    

\nwc{\cc}{\overline}   

\nwc{\IA}{\mathbb{A}} 
\nwc{\IB}{\mathbb{B}} 
\nwc{\IC}{\mathbb{C}} 
\nwc{\ID}{\mathbb{D}} 
\nwc{\IE}{\mathbb{E}} 
\nwc{\IF}{\mathbb{F}} 
\nwc{\IG}{\mathbb{G}} 
\nwc{\IH}{\mathbb{H}} 
\nwc{\IN}{\mathbb{N}} 
\nwc{\IP}{\mathbb{P}} 
\nwc{\IQ}{\mathbb{Q}} 
\nwc{\IR}{\mathbb{R}} 
\nwc{\IS}{\mathbb{S}} 
\nwc{\IT}{\mathbb{T}} 
\nwc{\IZ}{\mathbb{Z}} 

\nwc{\ba}{\blds{a}}
\nwc{\bb}{\blds{b}}
\nwc{\bc}{\blds{c}}
\nwc{\bd}{\blds{d}}
\nwc{\be}{\blds{e}}
\nwc{\fb}{\blds{f}} 
\nwc{\bg}{\blds{g}}
\nwc{\bh}{\blds{h}}
\nwc{\bi}{\blds{i}}
\nwc{\bj}{\blds{j}}
\nwc{\bk}{\blds{k}}
\nwc{\bl}{\blds{l}}
\nwc{\bm}{\blds{m}}
\nwc{\bn}{\blds{n}}
\nwc{\bo}{\blds{o}}
\nwc{\bp}{\blds{p}}
\nwc{\bq}{\blds{q}}
\nwc{\br}{\blds{r}}
\nwc{\bs}{\blds{s}}
\nwc{\bt}{\blds{t}}
\nwc{\bu}{\blds{u}}
\nwc{\bv}{\blds{v}}
\nwc{\bw}{\blds{w}}
\nwc{\bx}{\blds{x}}
\nwc{\by}{\blds{y}}
\nwc{\bz}{\blds{z}}

\nwc{\bA}{\blds{A}}
\nwc{\bB}{\blds{B}}
\nwc{\bC}{\blds{C}}
\nwc{\bD}{\blds{D}}
\nwc{\bE}{\blds{E}}
\nwc{\bF}{\blds{F}}
\nwc{\bG}{\blds{G}}
\nwc{\bH}{\blds{H}}
\nwc{\bI}{\blds{I}}
\nwc{\bJ}{\blds{J}}
\nwc{\bK}{\blds{K}}
\nwc{\bL}{\blds{L}}
\nwc{\bM}{\blds{M}}
\nwc{\bN}{\blds{N}}
\nwc{\bO}{\blds{O}}
\nwc{\bP}{\blds{P}}
\nwc{\bQ}{\blds{Q}}
\nwc{\bR}{\blds{R}}
\nwc{\bS}{\blds{S}}
\nwc{\bT}{\blds{T}}
\nwc{\bU}{\blds{U}}
\nwc{\bV}{\blds{V}}
\nwc{\bW}{\blds{W}}
\nwc{\bX}{\blds{X}}
\nwc{\bY}{\blds{Y}}
\nwc{\bZ}{\blds{Z}}

\nwc{\va}{{\bf a}}
\nwc{\vb}{{\bf b}}
\nwc{\vc}{{\bf c}}
\nwc{\vd}{{\bf d}}
\nwc{\ve}{{\bf e}}
\nwc{\vf}{{\bf f}}
\nwc{\vg}{{\bf g}}
\nwc{\vh}{{\bf h}}
\nwc{\vi}{{\bf i}}
\nwc{\vj}{{\bf j}}
\nwc{\vk}{{\bf k}}
\nwc{\vl}{{\bf l}}
\nwc{\vm}{{\bf m}}
\nwc{\vn}{{\bf n}}
\nwc{\vo}{{\it o}}
\nwc{\vp}{{\bf p}}
\nwc{\vq}{{\bf q}}
\nwc{\vr}{{\bf r}}
\nwc{\vs}{{\bf s}}
\nwc{\vt}{{\bf t}}
\nwc{\vu}{{\bf u}}
\nwc{\vv}{{\bf v}}
\nwc{\vw}{{\bf w}}
\nwc{\vx}{{\bf x}}
\nwc{\vy}{{\bf y}}
\nwc{\vz}{{\bf z}}

\nwc{\vA}{{\bf A}}
\nwc{\vB}{{\bf B}}
\nwc{\vC}{{\bf C}}
\nwc{\vD}{{\bf D}}
\nwc{\vE}{{\bf E}}
\nwc{\vF}{{\bf F}}
\nwc{\vG}{{\bf G}}
\nwc{\vH}{{\bf H}}
\nwc{\vI}{{\bf I}}
\nwc{\vJ}{{\bf J}}
\nwc{\vK}{{\bf K}}
\nwc{\vL}{{\bf L}}
\nwc{\vM}{{\bf M}}
\nwc{\vN}{{\bf N}}
\nwc{\vO}{{\it O}}
\nwc{\vP}{{\bf P}}
\nwc{\vQ}{{\bf Q}}
\nwc{\vR}{{\bf R}}
\nwc{\vS}{{\bf S}}
\nwc{\vT}{{\bf T}}
\nwc{\vU}{{\bf U}}
\nwc{\vV}{{\bf V}}
\nwc{\vW}{{\bf W}}
\nwc{\vX}{{\bf X}}
\nwc{\vY}{{\bf Y}}
\nwc{\vZ}{{\bf Z}}

\nwc{\cA}{\ml{A}}
\nwc{\cB}{\ml{B}}
\nwc{\cC}{\ml{C}}
\nwc{\cD}{\ml{D}}
\nwc{\cE}{\ml{E}}
\nwc{\cF}{\ml{F}}
\nwc{\cG}{\ml{G}}
\nwc{\cH}{\ml{H}}
\nwc{\cI}{\ml{I}}
\nwc{\cJ}{\ml{J}}
\nwc{\cK}{\ml{K}}
\nwc{\cL}{\ml{L}}
\nwc{\cM}{\ml{M}}
\nwc{\cN}{\ml{N}}
\nwc{\cO}{\ml{O}}
\nwc{\cP}{\ml{P}}
\nwc{\cQ}{\ml{Q}}
\nwc{\cR}{\ml{R}}
\nwc{\cS}{\ml{S}}
\nwc{\cT}{\ml{T}}
\nwc{\cU}{\ml{U}}
\nwc{\cV}{\ml{V}}
\nwc{\cW}{\ml{W}}
\nwc{\cX}{\ml{X}}
\nwc{\cY}{\ml{Y}}
\nwc{\cZ}{\ml{Z}}

\nwc{\nn}{\nonumber}
\nwc{\veta}{\mbox{\boldmath{$\eta$}}}
\nwc{\vhatk}{\hat{\bk}}
\nwc{\commentout}[1]{}
\nwc{\pdr}[2]{\frac{\partial{#1}}{\partial{#2}}}
\nwc{\odr}[2]{\frac{d{#1}}{d{#2}}}
\nwc{\summ}[2]{\sum_{#1=#2}^{\infty}}
\nwc{\summm}[3]{\sum_{#1=#2}^{#3}}
\nwc{\pdd}[1]{\frac{\partial}{\partial{#1}}}
\nwc{\odd}[1]{\frac{d}{d{#1}}}
\nwc{\m}{\mbox}
\nwc{\re}{\hbox{Re}}
\nwc{\mt}{\bar{t}}
\nwc{\ou}{\bar u}
\nwc{\noi}{\noindent}
\nwc{\non}{\nonumber}
\nwc{\half}{\frac{1}{2}}
\nwc{\third}{\frac{1}{3}}
\nwc{\pa}{\partial}
\nwc{\uP}{{\em \bf Proof: }}
\nwc{\uT}{\underline{Theorem:}}
\nwc{\epto}{{\xrightarrow{\ep\to 0}}}
\nwc{\nto}{{\xrightarrow{n\to \infty}}}

\begin{document}

\title{Dissipation time and decay of correlations}

\author{Albert Fannjiang\dag, 
St\'{e}phane Nonnenmacher\ddag\ and
Lech Wo{\l}owski\dag}
\thanks{\dag\ Department of Mathematics, 
University of California at Davis,
Davis, CA 95616, USA
(email: {\tt fannjian@math.ucdavis.edu}, {\tt wolowski@math.ucdavis.edu}).
The research of AF is supported in part by the grant from U.S. National
Science Foundation, DMS-9971322 and UC Davis Chancellor's Fellowship}
\thanks{\ddag\ Service de Physique Th\'eorique,
CEA/DSM/PhT (Unit\'e de recherche associ\'ee au CNRS) CEA/Saclay 91191
Gif-sur-Yvette c\'edex, France ({\tt nonnen@spht.saclay.cea.fr})}

\begin{abstract}
We consider the effect of noise on the dynamics generated by volume-preserving maps 
on a $d$-dimensional 
torus. The quantity we use to measure the irreversibility
of the dynamics is the {\sl dissipation time}.  We focus on the asymptotic behaviour of this
time in the limit of small noise. We derive universal lower
and upper bounds for the dissipation time in terms of various properties of the map and
its associated propagators: spectral properties, 
local expansivity, and global
mixing properties. We show that the dissipation is slow
for a general class of non-weakly-mixing maps; on the opposite, it is fast for
a large class of exponentially mixing systems which include
uniformly expanding maps and Anosov diffeomorphisms. 
\end{abstract}

\maketitle


\section{Introduction}

The origin of irreversibility in dynamical systems is often modeled by
small stochastic perturbations of the otherwise reversible dynamics. These
perturbations may be attributed to uncontrolled interactions with the
`environment', or with internal variables neglected in the equations. 
In experimental or numerical investigations, stochasticity or `noise' is introduced respectively
by finite precision
of the preparation and measurement apparatus, and by rounding-off errors
due to finite computer precision. 

One may take these stochastic perturbations explicitly into account by adding a term of
Langevin type in the evolution equations, or equivalently introducing some
diffusion term in the Fokker-Planck equation. In the present article
we choose to deal with discrete-time dynamics on some compact phase space, namely the 
$d$-dimensional torus. All the maps we will study preserve the Lebesgue measure, which is
therefore the ``natural'' invariant measure associated with the dynamics. 
The noise is implemented through a convolution operator, the kernel
of which has a (small) width $\ep>0$. For a given map $F$, two types of stochastic
perturbations will be considered: the noise operator may be applied at each step
of evolution, resulting in a ``noisy evolution''; on the opposite, we may choose to
introduce some stochasticity only at the initial (preparation) and final (measurement)
steps, resulting in a ``coarse-grained evolution''. 

In general, the influence of the noise on the long-time evolution of the system depends on the
dynamical properties of the map. Typically, the effect is stronger if the map
is `chaotic' than if it is `regular' \cite{Kif,B,N}. We will put this statement
into quantitative estimates.
To this end, we will characterize the
effect of noise on a map $F$ through a single observable, the {\sl dissipation time}. 
By dissipation, we mean the damping (through the noisy or coarse-grained evolution) of the
density fluctuations, while the total probability remains constant. Mathematically, this
damping is expressed by the decay of the $L^2$-norm of the density (or equivalently, the
$L^2$-norm of the density fluctuations) with respect the the Lebesgue measure. 
The decrease of this norm due to dissipation 
may be interpreted as an ``entropy increase'' of the system, up to the state of maximal entropy,
which is the uniform density
(here the word ``entropy'' denotes the Boltzmann entropy of the
probability density, unrelated with the topological or Kolmogorov-Sinai entropies associated
with the map).
We define the dissipation time as the time needed for the evolution
to bring fluctuations under a fixed threshold (that is, reduce the $L^2$ norm
by a fixed factor) \cite{FW}. We are especially interested in the behaviour of this time in the limit of 
small noise. We will show that this behaviour drastically depends on the dynamics generated
by the noiseless map $F$: in short, the dissipation will be ``fast'' for a chaotic dynamics, as opposed
to ``slow'' for a regular dynamics, the difference being ever more striking as the level of noise decreases.
Our main results exhibit this opposition by means of quantitative estimates on the asymptotics of the 
dissipation time.

Since the dissipation time is defined in terms of an $L^2$ norm, it is naturally related to
the spectral properties of the propagator associated with the map, 
acting on the space of square-integrable functions. 
In Section~\ref{SP}, we analyze the links between, on one side, the spectrum and pseudospectrum 
of the noisy or noiseless propagators, and on the other side, their dissipation time. The relevance of 
the pseudospectrum for time evolution problems has been recently put in evidence in the context of
non-unitary continuous-time dynamical systems \cite{Dav}.
These spectral relationships will be mainly used to analyze the case of ``regular'', precisely 
non-weakly-mixing maps. On the opposite, they are of little help for more chaotic ones.

In the following Sections, we connect the dissipation time to more `dynamical' properties of
the map $F$ (under some smoothness assumptions on $F$). We first obtain lower bounds for the 
dissipation time from the {\sl local expansion} properties of the map: a weak, or absence of
local expansion will imply slow dissipation (Section~\ref{LB}). Next, we obtain upper bounds
for the dissipation time using information on mixing properties of the map (that is, the
time decay of correlations between observables). Strong (e.g. exponential) mixing is a typical characteristics of 
chaotic behaviour and can be easily measured in numerical simulations \cite{Bal-decay}. The main
conclusion is that exponential or stronger mixing implies fast dissipation, for both the noisy and
the coarse-grained evolution.

In Section~\ref{examples}, we describe several families of volume-preserving maps on the torus, for which the
results obtained above yield useful information concerning the dissipation. The main two
families are, on the one hand, expanding (noninvertible) maps, on the other hand, Anosov (or at least
partially hyperbolic) diffeomorphisms of the torus. All these maps are at least exponentially mixing, 
so the dissipation is ``fast''. In case of Anosov maps, we collect our results in 
Theorem~\ref{TAnosov}. As a matter of illustration, we also analyze in detail examples of {\sl linear}
mixing maps, for which exact asymptotics of the dissipation times can be obtained, and therefore give
an idea of the sharpness of the bounds obtained before (the results concerning
 linear automorphisms had been obtained in \cite{FW}).

From these developments one notices that dissipation time
provides a robust characteristics of the chaoticity of a given dynamical system. As opposed to e.g. the 
decay rate of dynamical correlation, dissipation time has the same asymptotics (up to a constant factor) 
among Anosov diffeomorphisms while the decorrelation may be exponential
(generic Anosov case) or super-exponential depending on particular map (cf. results regarding 
hyperbolic toral automorphisms \cite{B} and Corollary \ref{supexp}). 


\section{Setup and notation}
\label{setup}

\subsection{Evolution operators}
\label{evop}
Let $(\IT^{d},\EuScript{B}(\IT^{d}),m)$ denote the $d$-dimensional torus, 
equipped with 
its $\sigma-$field of Borel sets and the Lebesgue measure
$m$.  Let $F:\IT^{d}\to \IT^{d}$ be a map on
the torus preserving the Lebesgue measure: for any set $B\in
\EuScript{B}(\IT^{d})$ we have $m(F^{-1}(B))=m(B)$. In general, $F$ is 
not supposed to be invertible. In the following we call such a map
`volume preserving' with implicit reference to the Lebesgue measure.

The map $F$ generates a discrete time dynamics on $\IT^{d}$,
which in terms of  pathwise description can be represented
by the forward trajectory 
$\{F^{n}(\bx_{0})$, $n\in \IN\}$ of any initial point (particle) $\bx_{0}\in \IT^{d}$.
However, instead of looking at the evolution 
of a single particle, one can consider the statistical description of
the dynamics, that is the evolution of a density (more generally 
a measure) describing the initial statistical configuration of the system. 

Let $\ml{M}(\IT^{d})$ denote the set of all Borel measures on $\IT^{d}$. 
For any $\mu\in\ml{M}(\IT^{d})$ and $f\in C^{0}(\IT^{d})$ we write
\bean
\mu(f)=\int_{\IT^{d}}f(\bx)d\mu(\bx).
\eean
The map $F$ induces a map $F^{*}$ on $\ml{M}(\IT^{d})$ given by
\bean
 (F^{*}\mu)(f)=\mu(f \circ F), \qquad \text{for all} \; f\in C^{0}(\IT^{d}).
\eean 
This map can also be defined as follows: 
\bean
 (F^{*}\mu)(B)=\mu(F^{-1}(B)), \qquad \text{for all} \; B\in \EuScript{B}(\IT^{d}).
\eean
In particular if $\mu=\del_{\bx_0}$ then $F^{*}(\mu)=\del_{F(\bx_0)}$ 
and one recovers the pathwise description. 

If $\mu$ is absolutely continuous w.r.t. $m$, then $F^{*}(\mu)$ preserves
this property (since the measure-preserving map $F$ is nonsingular w.r.t. $m$,
see \cite[p.42]{LM}). 
The corresponding densities $g=\frac{d\mu}{dm}\in L^1(\IT^d)$ 
are transformed by the Frobenius-Perron
or transfer operator $P_{F}$ \cite{B}:
\bean
P_F\left(\frac{d\mu}{dm}\right)=\frac{d(F^{*}\mu)}{dm}.
\eean

If the map $F$ is invertible, $P_F$ is given explicitly by:
\bean
\label{exF-P}
(P_{F}g)(\bx)= (g\circ F^{-1})(\bx)\frac{dF^{*}m}{dm}(\bx)
=g\circ F^{-1}(\bx).
\eean
If the map $F$ is differentiable, and the preimage
set of $\bx$ is finite for all $\bx$, the Perron-Frobenius operator is
given by
\bean
\label{0exF-P}
(P_{F}g)(\bx)=\sum_{\by|F(\by)=\bx}\frac{g(\by)}{|J_{F}(\by)|},
\eean
where $J_{F}(\by)$ is the Jacobian of $F$ at $\by$. 


On the other hand one can consider the dual
of the Frobenius-Perron operator, called
the Koopman operator, which governs the evolution of
observables
$f\in L^{\infty}(\IT^{d})$  instead of that of densities
$g\in L^1$. The   Koopman
operator 
$U_{F}$ is defined as
\bea
\label{exK}
U_{F}f=f\circ F.
\eea

Due to the nonseparability of the Banach space $L^{\infty}(\IT^{d})$, it
is often more convenient to consider its closure in some weaker $L^p$ norm, 
which yields larger (but separable) spaces of observables  
$L^{p}(\IT^{d})$. In this paper we will be mainly concerned with the space $L^2(\IT^d)$ and
its codimension-$1$ subspace of zero-mean functions $L^2_0(\IT^d)=\{f\in L^2(\IT^2) : m(f)=0\}$.
This subspace is obviously invariant under $U_F$ and $P_F$, due to the assumption
$F^*m=m$.
Throughout the paper, $\|\cdot\|$ will always refer to the $L^2$-norm (and corresponding
operator norm) on $L^2_0(\IT^d)$ (any other norm will carry an explicit subscript).

For any measure-preserving map $F$, the operator $U_F$ is isometric on $L^2(\IT^d)$
and $L^2_0(\IT^d)$. When $F$ is invertible, $U_F$ is unitary on these spaces, 
and satisfies $U_F=P_F^{-1}=P_{F^{-1}}$.


\subsubsection{Additional notation}
Although the operators introduced in previous sections were mostly defined on the space $L^2_{0}(\IT^d)$, 
it will be useful to consider other function spaces, which we define now in some detail.
For any $m\in\IN$, we denote by $C^m(\IT^d)$ the space of $m$-times continuously
differentiable functions, with the norm
\bean
\|f\|_{C^{m}}=\sum_{|\bmal|_{1}\leq m}\|D^{\bmal}f\|_{\infty}
\eean
(we use the norm $|\bmal|_1=\alpha_1+\ldots+\alpha_d$ for the multiindex $\bmal\in\IN^d$).
For any $s=m+\eta$ with $m=[s]\in\IN$, $\eta\in (0,1)$, let $C^{s}(\IT^d)$ denote the space of
$C^m$ functions for which the $m$-derivatives are $\eta$-H\"{o}lder continuous; this
space is equipped with the norm
\bean
\|f\|_{C^{s}}=\|f\|_{C^{m}}+
\sum_{|\bmal|_{1}= m}\sup_{\bx\not=\by}\frac{|D^{\bmal} f(\bx)-D^{\bmal} f(\by)|}{|\bx-\by|^\eta}
\eean 

The Fourier transforms of functions $g\in L^1(\IR^d)$ and $f\in L^1(\IT^d)$ 
are defined as follows:
\bea
\label{FT}
\forall \bxi\in \IR^d,\quad
\hat{g}(\bxi)&=& \int_{\IR^{d}}g(\bx)e^{-2\pi i \bx \cdot \bxi} d\bx, \\
\forall \bk\in \IZ^d,\quad
\hat{f}(\bk)&=& \int_{\IT^{d}}f(\bx)e^{-2\pi i \bx \cdot \bk} d\bx=\la\be_{\bk},f \ra.
\eea
Above we used the Fourier modes on the torus 
$\be_{\bk}(\bx):=e^{2 \pi i \bx\cdot\bk}$.
For any $s\geq 0$, we denote by $H^s(\IT^d)$ and $H^s(\IR^d)$ the Sobolev spaces of $s$-times
weakly differentiable $L^2$-functions equipped with the norms $\|\cdot\|_{H^s}$ defined
respectively by 
\bean
\|g\|^2_{H^s(\IR^d)}&=&\int_{\bxi\in \IR^d}(1+|\bxi|^2)^s |\hat{g}(\bxi)|^2d\bxi,\\
\|f\|^2_{H^s(\IT^d)}&=&\sum_{\bk\in \IZ^d}(1+|\bk|^2)^s |\hat{f}(\bk)|^2.
\eean
Finally, for any of these spaces, adding the subscript $0$ will mean that we consider the 
($U_F$-invariant) subspace
of functions with zero average, e.g. $C^j_0(\IT^d)=\{f\in C^j(\IT^d),\ m(f)=0\}$.


\subsection{Noise operator}
To construct the noise operator we first define the {\em noise generating density}
i.e. an arbitrary probability density function $g\in L^1(\IR^d)$ with even parity w.r.t. the origin:
$g(\bx)=g(-\bx)$.
The \emph{noise width}
(or \emph{noise level}) will be given by a single nonnegative parameter, which we call 
$\ep$. To each $\ep > 0$
corresponds the noise kernel on $\IR^d$:
\bean
g_{\ep}(\bx)=\frac{1}{\ep^{d}}g\left(\frac{\bx}{\ep}\right),
\eean
with the convention that $g_0=\del_0$. 
The noise kernel on the torus is obtained by periodizing $g_{\ep}$, yielding the
periodic kernel
\bea
\label{ank}
\tilde{g}_{\ep}(\bx)=\sum_{\bn\in \IZ^d} g_{\ep}(\bx+\bn).
\eea
We remark that the Fourier transform of $\tilde{g}_{\ep}$ is related to that of $g$ by
the identities $\hat{\tilde{g}}_{\ep}(\bk)=\hat{g}_{\ep}(\bk)=\hat{g}(\ep\bk)$.

The noise operator $G_{\ep}$ is defined on any function $f\in L^2_0(\IT^d)$ 
as the convolution:
\bean
G_{\ep}f &=& \tilde{g}_{\ep}*f.
\eean
As a convolution operator defined by an $L^1$ density, $G_{\ep}$ is compact on $L^2_0(\IT^d)$ 
(if $g$ is square-integrable, $G_\ep$ is Hilbert-Schmidt). 
The Fourier modes $\{\be_{\bk},\ \bk\in\IZ^d\setminus \{0\}\}$ form an orthonormal basis of
eigenvectors of $G_\ep$, yielding the following 
spectral decomposition:
\bea
\label{specdecompo}
\forall f\in L_0^2(\IT^d),\quad G_{\ep}f= \sum_{0\not=\bk\in \IZ^d} \hat{g}(\ep\bk) 
\la\be_{\bk},f \ra \,\be_{\bk}.
\eea
This formula shows that the eigenvalue associated with $\be_{\bk}$ is $\hat{g}(\ep\bk)$.
Since $g$ is a symmetric function, this eigenvalue is real, so that $G_\ep$ 
is a self-adjoint operator. Its spectral radius $r_{sp}(G_{\ep})$ is therefore 
given by
\bea
\label{Gspr}
r_{sp}(G_{\ep})=\|G_{\ep}\|=\sup_{0\not=\bk\in \IZ^d}|\hat{g}(\ep\bk)|.
\eea
Since the density $g$ is positive, $\hat g$ attains its maximum 
$\hat g(0)=1$ at the
origin and nowhere else. Besides, because $g\in L^1(\IR^d)$, 
$\hat g$ is a continuous function vanishing at infinity.  As a result, for small enough $\eps>0$, 
the supremum on the RHS of (\ref{Gspr}) is reached at some point $\ep\bk$ close to the origin, and 
this maximum is strictly smaller than $1$.
This shows that the operator $G_{\ep}$ is strictly contracting on $L_{0}^{2}(\IT^{d})$:
\bea
\label{strcon}
\forall \ep>0,\quad \|G_{\ep}\|=r_{sp}(G_{\ep})<1.
\eea
In the next section we study this noise operator more precisely, starting from appropriate
assumptions on the noise generating density.


\subsection{Noise kernel estimates}

In this subsection we present some estimates regarding the noise operator, which
will be used throughout the paper to estimate the dissipation time.

We will be interested in the behaviour of the system in the limit 
of small noise level, that is the limit $\ep\to 0$.
It will hence be useful to introduce the following asymptotic notation.
Given two variables $a_{\ep}$, $b_{\ep}$ depending on $\ep>0$,
we write 
\bea
a_{\ep} &\lesssim& b_{\ep} \text{ if } \limsup_{\ep
\rightarrow
0}\frac{a_{\ep}}{b_{\ep}} < \infty, \\
a_{\ep} &\approx& b_{\ep} \text{ if } \lim_{\ep \rightarrow 0}\frac{a_{\ep}}{ b_{\ep}}=1,\\
a_{\ep} &\sim& b_{\ep} \text{ if }a_{\ep} \lesssim b_{\ep}\text{ and }b_{\ep} \lesssim a_{\ep}.
 \label{asymnot}
\eea

In order to obtain interesting estimates on the noise operator $G_\ep$, it will be necessary to 
impose some additional conditions on its generating density $g$, regarding e.g. its rate of
decay at infinity, or the behaviour of its Fourier transform near the
origin.

The weakest condition considered in this paper is the existence of some
positive moment of $g$, by which we mean that for some $\alpha \in (0,2]$,
\bea
\label{Moment}
M_{\alpha}=\int_{\IR^d}|\bx|^{\alpha}g(\bx)d\bx < \infty
\eea
(we take the length $|\bx|=(x_1^2+\ldots+x_d^2)^{1/2}$ on $\IR^d$). 
This condition implies the following properties of the Fourier transform $\hat g$ (proved in
Appendix~\ref{app1}):

\begin{lem}
\label{Mom-Four}
For any $\alpha\in (0,2]$ there exists a universal constant $C_\alpha$ such that, 
if a normalized density $g$ satisfies (\ref{Moment}), then the following inequalities hold:
\bea\label{Moment-Fourier}
\forall \bxi\in\IR^d,\quad 0\leq 1-\hat g(\bxi)\leq C_\alpha M_\alpha |\bxi|^{\alpha}.
\eea
If (\ref{Moment}) holds with $\alpha=2$, we have the more precise behaviour:
$$
1-\hat g(\bxi)\sim |\bxi|^2\quad\mbox{in the limit}\ \bxi\to 0.
$$
\end{lem}

In the case $\alpha<2$, we will sometimes assume a stronger property than \eqref{Moment-Fourier}, namely that
\bea
\label{noise1}
 1-\hat{g}(\bxi) \sim |\bxi|^{\alpha}\quad \mbox{in the limit }\bxi\to 0.
\eea
Notice that this behaviour 
implies a uniform bound $1-\hat g(\bxi)\leq C|\bxi|^\gamma$ for any $\gamma\leq\alpha$ and $C$ independent of
$\gamma$.

Typical examples of noise kernels satisfying (\ref{noise1})
include the Gaussian kernel and more general symmetric $\alpha-$stable kernels
\cite[p.152]{Stroock} defined for $\alpha\in(0,2]$:
\bea
\label{Akernel}
g_{\ep,\alpha}(\bx):=\sum_{\bk \in \IZ^{d}}e^{-(\bQ(\ep \bk))^{\alpha/2}}\be_{\bk}(\bx),
\eea
where $\bQ$ denotes an arbitrary positive definite quadratic form. For the values
of $\alpha$ indicated, the function $g_{\ep,\alpha}(x)$ is positive on $\IR^d$.

\medskip

In view of Eq.\eqref{Gspr}, the properties \eqref{Moment} or \eqref{noise1}
directly constrain the rate at which  $G_{\ep}$ 
contracts on $L_{0}^{2}(\IT^{d})$. For instance, \eqref{noise1} implies that
in the limit $\ep\to 0$,
\bea
\label{contr}
1-\|G_{\ep}\| \sim \ep^{\alpha}.
\eea

The following proposition describes the effect of the noise on various types of
observables, in the limit of small noise level.
The proofs are given in Appendix~\ref{app2}. 

\begin{prop}
\label{nke}
i) For any noise generating density $g\in L^1(\IR^d)$ and any observable $f\in L^2_0(\IT^d)$, one has
\begin{equation}\label{gen}
\|G_{\ep}f-f\| \epto 0.
\end{equation}
To obtain information on the speed of convergence, we need to impose constraints on
both the noise kernel and the observable.

ii) If for some $\alpha\in(0,2]$ the kernel $g$ satisfies (\ref{Moment}) or (\ref{noise1}),
then for any $\gamma>0$ there exists a constant $C>0$ such that for any observable
$f\in H^{\gamma}(\IT^d)$,
\bea
\label{H1est}
\|G_{\ep}f-f\| \leq  C \ep^{\gamma \wedge \alpha} \|f\|_{H^{\gamma \wedge \alpha}},
\eea
where $ \gamma\wedge \alpha:=\min\{\gamma,\alpha\}$.
If $f\in C^{1}(\IT^d)$, the above upper bound can be replaced by
\bea
\label{C1est}
\|G_{\ep}f-f\| \leq C \ep^{1\wedge \alpha} \|\nabla f\| 
\leq C \ep^{1\wedge\alpha }\|\nabla f\|_{\infty}.
\eea
\end{prop}

Using the noise operator, we are now in position to define the noisy (resp. the
coarse-grained) dynamics generated by a measure-preserving map $F$.


\subsection{Noisy evolution operator and dissipation time}
The noisy evolution through the map $F$ is constructed by successively applying
the Koopman operator $U_{F}$ and the noise operator $G_{\ep}$, therefore by
taking powers of the \emph{noisy propagator}
\bean
T_{\ep}=G_{\ep}U_{F}.
\eean
In general, the operator $T_{\ep}$ is not normal, but satisfies 
$r_{sp}(T_\ep)\leq \|T_\ep\|=\|G_\ep\|$. 

We will also consider a \emph{coarse-grained dynamics} 
defined by the application of the noise kernel only at
the beginning and the end of the evolution. Hence we define
the following family of operators:
\bean
\tilde{T}^{(n)}_{\ep}=G_{\ep}U^{n}_{F}G_{\ep},\quad n\in\IN.
\eean

In view of the contracting properties of $G_{\ep}$,
the inequalities $\|T_\ep^{n}\|\leq \|G_\ep\|^n$, $\|\tilde{T}_\ep^{(n)}\|\leq \|G_\ep\|^2$ 
imply that 
both noisy and coarse-grained operators are strictly contracting on $L_{0}^2(\IT^d)$.
Our aim is to characterize the speed of contraction of these two evolutions, that is the
behaviour of the norms $\|T_\ep^{n}\|$, $\|\tilde{T}_\ep^{(n)}\|$ in the joint
limits  $\ep\to 0$ and $n \to \infty$.
This characterization will be connected with dynamical properties of the map $F$.

There are
many ways to measure this speed of contraction. For instance, for fixed $\ep>0$, 
the long-time decay of $\|T_\ep^{n}\|$ may be super-exponential (in which
case $T_\ep$ is quasinilpotent) or exponential,  
governed e.g. by the largest eigenvalue of $T_\ep$. However, such exponential behaviour may
appear only after a transient time. 
In this paper we will characterize the noisy dynamics 
through a single, robust characteristics, namely the {\em dissipation time}.
In its general form the dissipation time ${\tau}_*$ 
is defined in terms of the norm $\|\cdot\|_{p,0}$
on the space $L^p_0(\IT^d)$, and an arbitrary threshold $\eta\in(0,1)$:
\bea
\label{tdiss2}
{\tau}_*^{p,\eta}(\ep):=\min\{n \in \IN:\|T_{\ep}^{n}\|_{p,0}< \eta\}, \quad 1\leq p\leq\infty.
\eea
We will be concerned with the behaviour of the dissipation time when the level of 
noise becomes small (as we prove in Proposition~\ref{nonfinite}, this time
diverges in this limit).
In \cite{FW} it was shown that this asymptotic behaviour is independent of 
the choice of $0<\eta<1$ 
and $1<p<\infty$. Therefore, in the present paper 
we will consider the computationally convenient choice
$p=2$ and $\eta=e^{-1}$. 
We will henceforth drop the superscripts, and the dependence of ${\tau}_*$ on $\ep$ 
will always be implicit:
\bea
\label{tdiss1}
{\tau}_*:={\tau}_*^{2,e^{-1}}(\ep)=\min\{n \in \IN:\|T_{\ep}^{n}\| < e^{-1}\}.
\eea

A similar dissipation time will be
defined for the coarse grained evolution:
\bea
\label{ttdiss1}
\tilde{\tau}_*:=\min\{n \in \IN:\|\tilde{T}_{\ep}^{(n)}\| < e^{-1}\}.
\eea

The dissipation time 
does not depend on whether the dynamics is applied to
densities (i.e. by the Frobenius-Perron operator) or to
observables (by the Koopman operator). 
Indeed, the norm of
an operator equals the norm of its adjoint \cite[p.195]{Y}, so that
\bean
\|\tilde{T}_{\ep}^{(n)}\|=
\|G_{\ep}U^{n}_{F}G_{\ep}\|=
\|(G_{\ep}U^{n}_{F}G_{\ep})^{*}\|
=\|G_{\ep}P^{n}_{F}G_{\ep}\|,
\eean
and similarly for the noisy operator $T_\ep$.
In particular, for invertible maps the dissipation time does not depend
on the direction of time.

\medskip

We will distinguish two qualitatively different asymptotic behaviours of dissipation time
in the limit $\ep\to 0$.
We say that the operator $T_{\ep}$ (or the map $F$ associated with it) respectively has
\begin{itemize}
\item[I)] {\em simple} or {\em power-law} dissipation time if there exists $\beta>0$ such that
\bean
{\tau}_* \sim 1/\ep^{\beta},
\eean
\item[II)] {\em fast} or {\em logarithmic} dissipation time if 
\bean
{\tau}_* \sim \ln(1/\ep).
\eean
\end{itemize}
We will also talk about {\em slow} dissipation time whenever there exists 
some $\beta>0$ s.t.
\bean
{\tau}_* \gtrsim 1/\ep^{\beta}.
\eean
In case of logarithmic dissipation time, the dissipation rate constant $R_*$, when it exists, is 
defined as
\bea
\label{Rdiss}
  R_*=\lim_{\ep \rightarrow 0} \frac{{\tau}_*}{\ln(1/\ep)}.
\eea    
A similar terminology will be applied to the coarse-grained dissipation time $\tilde{\tau}_{*}$.


\section{Spectral properties and dissipation time of non-weakly-mixing maps}
\label{SP}

In this section we investigate the connection between the dissipation time of the noisy propagator $T_{\ep}$ and 
its pseudospectrum together with some spectral properties of $U_{F}$ and $G_{\ep}$. All the operators considered
in this section are defined on $L^{2}_{0}(\IT^d)$. In the framework of continuous-time 
dynamics, connections have been obtained between, on one side, the
pseudospectrum of the (non-selfadjoint) generator $A$, and on the other side, the norm of the evolution operator 
$e^{tA}$ \cite{Dav}. We will obtain results of the same flavor, yet the proofs seem here easier
than in the case of continuous time. 


\subsection{Definitions and general bounds}
Let us start with the definition of the pseudospectrum of a bounded operator \cite{Var}.

\begin{defin}
\label{psp}
Let $T$ be a bounded linear operator on a Hilbert space $\mathcal{H}$ (we note $T\in \mathcal{L}(\mathcal{H})$).
For any $\del>0$, the $\del$-pseudospectrum 
of $T$ (denoted by $\sigma_{\del}(T)$) can be defined in the following three equivalent ways:
\bean
&(I)& \sigma_{\del}(T)=\{\lam \in \IC: \|(\lam-T)^{-1}\|\geq \del^{-1}\}, \\
&(II)& \sigma_{\del}(T)=\{\lam \in \IC: \exists v\in\mathcal{H}, \ \|v\|=1, \quad \|(T-\lam)v\|\leq \del \}, \\
&(III)& \sigma_{\del}(T)=\{\lam \in \IC: \exists B\in \mathcal{L}(\mathcal{H}),\ 
\|B\|\leq \del, \quad \lam \in \sigma(T+B) \}.
\eean
\end{defin}
We will apply these definitions to the operator $T_\ep$. 
For brevity, the resolvent of this operator will be denoted by 
$R_\ep(\lam)=(\lam - T_\ep)^{-1}$.
We call $S^{r}$ the circle $\{\lam \in \IC: |\lam|=r\}$ in the complex plane, and define the following 
{\em pseudospectrum distance function}:
\bean
d_\ep(r):=\inf \{ \del>0 :  \sigma_{\del}(T_\ep)\cap S^r \not = \emptyset\}.
\eean
From the definition (I) of the pseudospectrum, one easily shows that
this distance is also given by
\bea
\label{psdist}
d^{-1}_\ep(r)=\sup_{|\lam|=r} \|R_\ep(\lam)\|.
\eea

We first establish general (abstract) bounds for the dissipation time in terms of
the spectral properties of $G_{\ep}$ and $T_{\ep}$. In a second step, we relate these properties
to dynamical properties of the underlying map $F$.

\begin{thm}
\label{gb}
For any isometric operator $U$ on $L^2_0(\IT^d)$ and noise operator $G_\ep$,
the dissipation time of the noisy evolution operator $T_\ep=G_{\ep}U$ satisfies the following estimates:
\begin{align}
\frac{1-e^{-1}}{d_\ep(1)} \leq  {\tau}_* &\leq
\frac{1}{\left|\ln(\|G_\ep\|)\right|}+1,\label{firstline}\\
 {\tau}_* &\leq\inf_{r_{sp}(T_\ep)<r<1}\frac{1}{|\ln (r)|} \ln\left(\frac{e}{d_{\ep}(r)}\right).
\label{second-upper}
\end{align}
\end{thm}

We notice that the first upper bound does not depend on $U$ at all, but only on the noise.
Using the estimate \eqref{contr}, we obtain the following obvious corollary:

\begin{cor}
\label{noise1->upper}
If the noise generating density satisfies the
estimate \eqref{noise1} for some $\alpha\in (0,2]$, then for any measure-preserving map $F$ the
noisy dissipation time is bounded from above by:
\bean
{\tau}_* \lesssim  \ep^{-\alpha}.
\eean
\end{cor}

\textbf{Proof of the Theorem.} 

1. \underline{Lower bound}

We use the following series expansion of the resolvent \cite[p.211]{Y}
valid for any $|\lam|> r_{sp}(T_\ep)$:
\bea
\label{NSR}
R_\ep(\lam)=\sum_{n=0}^{\infty}\lam^{-n-1}T_\ep^{n}.
\eea 
Considering that $r_{sp}(T_{\ep})\leq \|G_{\ep}\|< 1$, we may take  $|\lam|=1$, and cut this sum
into two parts:
\bean
 R_{\ep}(\lam)=
\sum_{n=0}^{{\tau}_*-1}\lam^{-n-1}T_{\ep}^{n}+
\lam^{-{\tau}_*}T_{\ep}^{{\tau}_*}R_{\ep}(\lam).
\eean
Taking norms and applying the triangle inequality, we get
\bean
\|R_{\ep}(\lam)\| &\leq& \|\sum_{n=0}^{{\tau}_*-1}\lam^{-n-1}T_{\ep}^{n}\|
+|\lam|^{-{\tau}_*}\|T_{\ep}^{{\tau}_*}\|\|R_{\ep}(\lam)\|\\
&\leq& {\tau}_*+e^{-1}\|R_{\ep}(\lam)\|\\
\Longrightarrow \|R_{\ep}(\lam)\|(1-e^{-1})&\leq& {\tau}_*. 
\eean
Taking the supremum over $\lambda\in S_1$ yields the lower bound. 

\bigskip

2. \underline{Upper bounds}

To get both upper bounds, we use the following trivial lemma.

\begin{lem}\label{simple}
Assume that (for some value of $\ep$) the powers of $T_\ep$ satisfy
$$
\forall n\in\IN,\quad \|T_\ep^n\|\leq \Gamma(n),
$$
where the function $\Gamma(n)$ is strictly decreasing, and $\Gamma(n)\nto 0$.
Then the dissipation time is bounded from above by
$$
{\tau}_*\leq \Gamma^{(-1)}(e^{-1})+1,
$$
where $\Gamma^{(-1)}$ is the inverse function of $\Gamma$.
In particular, for the geometric decay
$\Gamma(n)=C r^n$ with $r\in(0,1)$,
$C \geq 1$, one obtains  ${\tau}_*\leq
\frac{\ln(eC)}{|\ln r|}+1$.
\end{lem}

The upper bound in Eq.~\eqref{firstline} comes from the obvious estimate
\bean
\|T_{\ep}^n\| \leq \|G_{\ep}\|^n,
\eean
on which we apply the lemma with $C=1$, $r=\|G_\ep\|$.

To prove the second upper bound, we use the representation of $T_\ep^{n}$ in terms of the resolvent:
\bean
T_\ep^{n}=\frac{1}{2\pi i}\int_{S^r}\lam^{n}R_{\ep}(\lam)d\lam
\eean
valid for any $r>r_{sp}(T_{\ep})$. Thus for all $r\in(r_{sp}(T_\ep),1)$, one has
\bean
\|T_\ep^{n}\|\leq \frac{1}{2\pi}\int_{S^{r}}|\lam|^{n}\|R_{\ep}(\lam)\||d\lam|
\leq \sup_{|\lam|=r}\|R_{\ep}(\lam)\| r^{n+1}=\frac{1}{d_\ep(r)}\,r^{n+1}.
\eean
We then apply Lemma \ref{simple} on the geometric
decay for any radius $r_{sp}(T_\ep)<r<1$, with $C=\frac{r}{d_\ep(r)}\geq 1$.

\hfill$\blacksquare$


\subsection{Consequences}
\label{conseq}
We now use Theorem~\ref{gb} in the case where $U=U_F$ is the Koopman operator for some measure-preserving
map $F$ on $\IT^d$ with some specific dynamical properties. 
We recall \cite{CoFoSi} that the map $F$ is ergodic (resp. weakly-mixing) iff $1$ is not an
eigenvalue of $U_{F}$ (resp. iff $U_{F}$ has no eigenvalue) on $L^2_0(\IT^d)$.

\begin{prop}\label{nonfinite}
For any measure-preserving map $F$ and any noise generating function $g$, the dissipation
time of $T_\ep$ diverges in the small-noise limit $\ep\to 0$. 
\end{prop}

\textbf{Proof.} We skip the subscript $F$ to alleviate notation. 
We only use the fact that $U=U_F$ is an isometry. 
We shall prove by induction the following strong convergence of operators
$$
\forall f\in L^2_0(\IT^d),\quad \forall n\in\IN,\qquad \|T_\ep^n f-U^n f\|\epto 0.
$$
From Proposition~\ref{nke}{\em i)}, this limit holds in the case $n=1$. 
Let us assume it holds at the rank $n-1$.
Then we write
$$
T_\ep^n f=UT_\ep^{n-1} f + (G_\ep -I)UT_\ep^{n-1} f.
$$
From the inductive hypothesis, $T_\ep^{n-1}f\epto U^{n-1}f$, so that the first term on the
RHS converges to $U^nf$. Applying Proposition~\ref{nke}{\em i)} to the
function $U^n f$, we see that the second term vanishes in the limit $\ep\to 0$.
From the isometry of $U$, we obtain that for any $n>0$, $\|T_\ep^n\|\epto 1$, 
so that ${\tau}_*\epto\infty$.\hfill$\blacksquare$

This ``non-finiteness'' of the noisy dissipation time allows us to prove another general result concerning
the pseudospectrum of $T_\ep$ (this corollary is proved in the Appendix~\ref{app3}):

\begin{cor}\label{disto0}
For any isometry $U$ and noise generating function $g$, one has
\bea
\label{depto0}
d_\ep(1)\epto 0.
\eea
This means that for any fixed $\delta>0$, the pseudospectrum $\sigma_\delta(T_\ep)$ will intersect
the unit circle for small enough $\ep$.
\end{cor}

In order to better control the growth of ${\tau}_*$, we need more 
precise information on the noise and the dynamics. In the present section, we restrict ourselves to
the dynamical property of weak-mixing.

\begin{cor}\label{weakmix}
Assume that the noise generating density $g$ satisfies the estimates  
\eqref{Moment} or \eqref{noise1} with exponent $\alpha\in (0,2]$.
If $F$ is not weakly-mixing and at least one eigenfunction of 
$U_{F}$ belongs to $H^{\gamma}(\IT^d)$ for some $\gamma>0$, 
then $T_{\ep}$ has slow dissipation time: 
\bean
\ep^{-(\alpha\wedge\gamma)} \lesssim {\tau}_*.
\eean 
\end{cor}
\textbf{Proof.} 
Let $h\in H^\gamma(\IT^d)$ be a normalized eigenfunction of $U_F$ with eigenvalue $\lam$. 
Applying Proposition~\ref{nke}{\em ii)}, we get
$$
\|(\lam -T_{\ep})h\|= \|(I-G_{\ep})h\|\leq K\ep^{\gamma\wedge\alpha}
$$
for some constant $K>0$ depending on $g$ and $h$. 
This implies that $\|R_\ep(\lam)\|\geq\frac{1}{K\ep^{\gamma\wedge\alpha}}$, therefore
taking the supremum over $|\lam|=1$ yields $d_\ep(1)^{-1}\geq\frac{1}{K\ep^{\gamma\wedge\alpha}}$.
The lower bound in Theorem~\ref{gb} then implies 
\bea
\label{low-gamma}
\frac{1-e^{-1}}{K\ep^{\gamma\wedge \alpha}} \leq  {\tau}_*.
\eea
\hfill$\blacksquare$

\begin{rem}
Recall that if $g$ satisfies \eqref{noise1} with exponent
$\alpha$, then  the dissipation time is also bounded from
above, as shown in Corollary~\ref{noise1->upper}. If one eigenfunction $h$ has
regularity $H^\gamma$ with $\gamma\geq\alpha$, then both Corollaries
imply that the dissipation is {\em simple}, of exponent
$\alpha$.
\end{rem}

\begin{rem}
The above results can be stated in more general form: $U_{F}$ does not
need to be a Koopman operator associated with a map $F$. The result
holds true for any isometric operator $U$ on $L^2_0$ with an eigenfunction of Sobolev regularity.
\end{rem}

The dependence of the lower bound in (\ref{low-gamma}) on $\gamma$ can be intuitively 
explained as follows. In case of non-weakly-mixing maps the eigenfunctions of $U_F$
are, in general, responsible for slowing down the dissipation. The rate of the dissipation 
is affected by the regularity of the smoothest eigenfunction. In principle, irregular functions
undergo faster dissipation giving rise to slower asymptotics of $\tau_{*}$.
It is not clear, however, whether the actual asymptotics of the dissipation time
will be slower than power law in case when all eigenfunctions of $U_F$ on $L^2_0(\IT^d)$ are
outside any space $H^{\gamma}(\IT^d)$ with $\gamma>0$.

In Corollary~\ref{disto0} we have shown that for any map $F$ and
arbitrary small $\del>0$, 
the pseudospectrum $\sigma_\del(T_\ep)$ 
intersects the unit circle for sufficiently small $\ep>0$.
If $F$ is not weakly-mixing, the spectral radius of $T_\ep$ (that is, the modulus of its largest eigenvalue) 
is believed to converge to $1$ when $\ep\to 0$, and the associated eigenstate $h_\ep$ should converge 
to a ``noiseless eigenstate'' $h$.
This ``spectral stability'' has been discussed for several cases in the continuous-time as well as
for discrete-time maps on $\IT^2$ \cite{Kif,N}.

On the opposite, if $F$ is an Anosov map on $\IT^2$ (see Section~\ref{Anosov}), 
the spectrum of $T_{\ep}$ does not approach the unit circle, but
stays away from it uniformly: $r_{sp}(T_\ep)$
is smaller than some $r_0<1$ for any $\ep>0$ \cite{BKL}. Simultaneously, $\|T_\ep\|\to 1$,
so we have here
a clear manifestation of the {\em nonnormality} of $T_\ep$ for such a map. In some
cases (see \cite{N} and the linear examples of Section~\ref{examples}), the operator $T_{\ep}$ is even 
quasinilpotent, meaning that $r_{sp}(T_{\ep})=0$ for all $\ep>0$. 
For such an Anosov map, the spectral radius of $T_\ep$ is therefore
``unstable'' or ``discontinuous'' in the limit $\ep\to 0$, while in the same limit 
the (radius of its) pseudospectrum
$\sigma_\del(T_\ep)$ (for $\delta>0$ fixed) is  ``stable''.


\section{Universal lower bounds for the dissipation time of $C^1$ maps}
\label{LB}
In this section we consider both noisy and coarse-grained evolutions.
We start our discussion with general properties of the coarse-grained dissipation time.

\begin{prop}\label{infinite} 
Let $F$ be a measure-preserving map. 

i) For arbitrary noise kernel, the coarse-grained dissipation time 
diverges as $\ep \rightarrow 0$.

ii) If $F$ is not weakly-mixing then $\tilde{\tau}_*=\infty$ for small enough $\ep>0$.
\end{prop}
\textbf{Proof.}

{\em i)} Similarly as in Proposition \ref{nonfinite} we have 
\bean
\|\tilde{T}_{\ep}^n f - U^{n}f\|=\|G_{\ep}U_{F}^n(G_{\ep}-I)f+(G_{\ep}-I)U^nf\|
\leq \|(G_{\ep}-I)f\|+\|(G_{\ep}-I)U^nf\|\rightarrow 0.
\eean

{\em ii)} Let $h\in L^2_0(\IT^d)$ be a normalized eigenfunction of $U_F$, then
\bean
\|\tilde{T}^{(n)}_{\ep}h\| &=& \|G_{\ep}U^n_F(h+(G_\ep -I)h)\|
\geq\|G_{\ep}h\|-\|G_{\ep}U^n_F(G_\ep -I)h\|\\
&\geq& 1-2\|(G_\ep -I)h\|.
\eean
Since the RHS above is independent of $n$, we see that $\|\tilde{T}^{(n)}_{\ep}\|$ is close to $1$ for
all times and sufficiently small $\ep>0$. \hfill$\blacksquare$

As opposed to the noisy case (see Prop.~\ref{noise1->upper}), the coarse-grained evolution through a
non-weakly-mixing map 
does not dissipate. Therefore there exists no general upper bound for $\tilde {\tau}_*$.

On the opposite, we will prove below a general {\em lower} bound for both coarse-grained and noisy evolutions,
valid for any measure-preserving map $F$ of regularity $C^1$.
We note that Propositions~\ref{nonfinite} and \ref{infinite}{\em i)} (which are valid
independently of any regularity assumption) do not provide an explicit lower bound.

First we introduce some notation. For any map $F\in C^1$, 
$DF(\bx)$ is the tangent map of $F$ at the point $\bx\in\IT^d$, mapping
a tangent vector at $\bx$ to a tangent vector at $F(\bx)$. Selecting the canonical (i.e. Cartesian) basis 
and metrics on $T(\IT^d)$, this map can be represented as a $d\times d$ matrix. 
The metrics naturally yields a norm $\bv\in T_{\bx}(\IT^d)\mapsto |v|$ on the tangent space, and
therefore a norm on this matrix: $|DF(\bx)|=\max_{|\bv|=1}|DF(\bx) \cdot \bv|$.
We are now in position to define the maximal expansion rate of $F$:
\bean
\mu_{F}=\limsup_{n\rightarrow \infty} \|DF^n\|_{\infty}^{1/n},\quad\mbox{where}\quad
\|DF^n\|_{\infty}=\sup_{\bx\in\IT^d} |(DF^n)(\bx)|.
\eean
Since $F$ preserves the Lebesgue measure, the Jacobian $J_F(\bx)$ satisfies
$|J_F(\bx)|\geq 1$ at all points. 
In the Cartesian basis, $J_F(\bx)=\det(DF(\bx))$, so
that we have $\|DF^n(\bx)\|\geq 1$ for all $\bx\in\IT^d$, $n\geq 0$. One can
actually prove the following:
\begin{rem}
Although $|(DF^n)(\bx)|$ and $\|DF\|_{\infty}$ may depend on the choice of the metrics, $\mu_F$ does not, 
and satisfies $1\leq \mu_F\leq \|DF\|_{\infty}$.
\end{rem}

From the definition of $\mu_{F}$, for any
$\mu > \mu_F$ there exists
a constant $A\geq 1$ such that
\bea
\label{expogrowth}
\forall n\in \IN,\qquad \|DF^n\|_{\infty} \leq A\mu^n.
\eea
In some cases one may take $\mu=\mu_F$ in the RHS.
In case $\mu_F=1$, $\|DF^n\|_{\infty}$ can sometimes grow as a power-law:
\bea
\label{power-law}
\|DF^n\|_{\infty} \leq A n^\beta, \qquad n\in \IN
\eea
for some $\beta>0$, or even be uniformly bounded by a constant ($\beta=0$). 

\medskip

The relationship between, on one side, the local expansion of the map $F$ and on the other side,
the dissipation time, can be intuitively understood as follows. A lack of 
expansion ($\|DF\|_\infty=1$) results in the transformation of ``soft'' or ``long-wavelength'' oscillations
into ``soft oscillations'', both being little affected by the
noise operator $G_\ep$. On the opposite, a locally strictly expansive map ($\|DF\|_\infty>1$) will
quickly transform soft oscillations into ``hard'' or ``short-wavelength'', 
the latter being much more damped by the noise.

The following theorem precisely measures this relationship, in terms of {\em lower bounds}
for the dissipation times.

\begin{thm}
\label{nln}
Let $F$ be a measure-preserving $C^1$ map on $\IT^d$, and assume that the 
noise generating density $g$ satisfies (\ref{Moment}) or
(\ref{noise1}) for some $\alpha\in (0,2]$.

i) If $\|DF\|_{\infty}> 1$, resp. $\mu_{F}>1$, then
there exist a constant $c$, resp. constants $\mu\geq\mu_F$ and $\tilde{c}$, such that for small enough $\ep$,
\bea
\label{nlb}
{\tau}_* \geq \frac{\alpha\wedge 1}{\ln(\|DF\|_\infty)} \ln(\ep^{-1})+c, \quad\mbox{resp.}\quad
\tilde{\tau}_* \geq \frac{\alpha \wedge 1}{\ln\mu} \ln(\ep^{-1})+\tilde{c}.
\eea
If $F$ is a $C^1$ diffeomorphism, then (\ref{nlb}) holds
with
$\|DF\|_{\infty}$ replaced by
$\|DF\|_{\infty}\wedge\|D(F^{-1})\|_{\infty}$, resp. with some
$\mu\geq\mu_{F}\wedge\mu_{F^{-1}}$.

ii) If $\|DF\|_{\infty}=1$ then $T_{\ep}$ has slow 
dissipation time, ${\tau}_*\gtrsim \ep^{-(\alpha\wedge 1)}$. 
If the noise kernel satisfies the condition (\ref{noise1}) for $\alpha\in (0,1]$, then 
the dissipation time is simple, $\tau_*\sim \ep^{-\alpha}$.

iii) If $\mu_{F}=1$ and $\|DF^n\|_\infty$ grows as a power-law as in Eq.~\eqref{power-law}
with $\beta>0$, then $\tilde{\tau}_* \gtrsim \ep^{-(\alpha\wedge 1)/\beta}$.
If $\|DF^n\|_\infty$ is uniformly bounded above by a constant, then $\tilde{\tau}_*=\infty$ 
for small enough $\eps$.
\end{thm}

\medskip

\begin{rem}
\label{class}
This theorem shows that classical systems on $\IT^d$ (i.e. $C^1$ diffeomorphisms) cannot have
a dissipation time growing slower than $C\ln(\ep^{-1})$. In view of the results 
for toral automorphisms (cf. Proposition~4), this lower bound
on the dissipation time is sharp and consistent with Kouchnirenko's upper
bound on the entropy of the classical systems, namely
all classical systems have a finite (possibly zero) Kolmogorov-Sinai
entropy (Theorem 12.35. in \cite{Arnold}, see also \cite{AM}, \cite{Kou}).
\end{rem}

\textbf{Proof of the Theorem.}

We will need the following trivial lemma (similar with Lemma~\ref{simple}).

\begin{lem}\label{simple2}
Assume that there exists some $\alpha>0$ and a strictly increasing function $\gamma(n)$, $\gamma(0)=0$ such that
\bea
\forall n\geq 1,\qquad \|T_\ep^n\|\geq 1-\ep^\alpha \gamma(n).
\eea
Then the dissipation time is bounded from below as:
\bea
{\tau}_*\geq \gamma^{(-1)}\left(\frac{1-e^{-1}}{\ep^\alpha}\right),
\eeq
where $\gamma^{(-1)}$ is the inverse function of $\gamma$.

The same statement  holds for  the coarse-grained version.
\end{lem}

Our task is therefore to bound $\|T_\ep^n\|$ (resp. $\|\tilde T_\ep^{(n)}\|$) from below.
A simple computation shows that for any $f\in C^0(\IT^d)$, 
$\|G_{\ep}f\|_{\infty}\leq \|f\|_{\infty}$. Since convolution commutes
with differentiation, for $f\in C^1$ we also have $\|\nabla(G_\ep f)\|_\infty\leq \|\nabla f\|_\infty$. 
We use this fact to estimate the gradient of $T_{\ep}f$:
\bean
\|\nabla (T_{\ep}f) \|_{\infty}&=&\|\nabla(G_{\ep}U_{F}f)\|_{\infty}\\
&\leq& \|\nabla(f\circ F) \|_{\infty}=\|(\nabla f)\circ F \cdot DF \|_{\infty}\\
&\leq& \|(\nabla f)\circ F\|_{\infty}\|DF\|_{\infty}= \|\nabla f\|_{\infty}\|DF\|_{\infty}.
\eean
Repeating the above procedure $m$ times, we get
\bea
\label{nab}
\|\nabla(T^{m}_{\ep}f) \|_{\infty}
\leq \|\nabla f\|_{\infty}\|DF\|_{\infty}^m,\quad 
\|\nabla(U_FT^{m}_{\ep}f) \|_{\infty}
\leq \|\nabla f\|_{\infty}\|DF\|_{\infty}^{m+1}.
\eea
We now choose some arbitrary $f\in C^1_0(\IT^d)$, with $\|f\|=1$.
We first apply the triangle inequality:
$$
\|T^n_{\ep}f\|= \|G_{\ep}U_{F}T^{n-1}_{\ep}f\|
\geq  \|U_{F} T^{n-1}_{\ep}f\| - \|(G_{\ep}-I)U_{F}T^{n-1}_{\ep}f\|.
$$
To estimate the second term on the RHS we use the bound \eqref{C1est} and
the estimate \eqref{nab} to obtain
$$
\|T^n_{\ep}f\|\geq \|T^{n-1}_{\ep}f\|-C\ep^{\alpha \wedge 1}\|\nabla f\|_{\infty}\|DF\|_{\infty}^n.
$$
Applying the same procedure iteratively to the first term on the RHS, we finally get (remember $\|f\|=1$):
\bea\label{est-noise}
\|T^n_\ep\|\geq\|T^n_{\ep}f\|\geq 1-C\ep^{\alpha \wedge 1} \|\nabla f\|_{\infty}\sum_{m=1}^n\|DF\|_{\infty}^m.
\eea
The computations in the case of the coarse-grained operator are even simpler:
\bea
\|\tilde{T}^{(n)}_{\ep}f\|&=&\|G_{\ep}U_{F}^n G_{\ep}f\|\non\\
&\geq&
1-C\ep^{\alpha\wedge 1} \|\nabla f\|_{\infty} - C\ep^{\alpha\wedge 1} \|\nabla(G_{\ep}f)\|_{\infty}\|DF^n\|_{\infty}
\non\\
&\geq&
1- 2C\ep^{\alpha\wedge 1} \|\nabla
f\|_{\infty}\|DF^n\|_{\infty}.\label{est-coarse}
\eea
Notice that from the assumptions on $f$, 
$\|\nabla f\|_\infty$ cannot be made arbitrary small, but is necessarily larger than some positive constant.
We choose some arbitrary function, say $f=\be_{\bk}$ with $\bk=(1,0)$ which satisfies 
$\|\nabla f\|_\infty=2\pi$.

\medskip

The estimate (\ref{est-noise}) has the form given in Lemma~\ref{simple2}. 
The growth of the function $\gamma(n)$ depends
on whether $\|DF\|_{\infty}$ is equal to or larger than $1$, 
which explains why the lower bounds are qualitatively different in the two cases. 
 
In case $\|DF\|_{\infty}$ is strictly larger than $1$, then the function $\gamma(n)$ 
grows like an exponential, therefore the lower bound is of the type \eqref{nlb}.
For the coarse-grained version, a growth of $\|DF\|_{\infty}$ of the type \eqref{expogrowth} 
yields the lower bound for $\tilde{\tau}_*$ in \eqref{nlb}.

In the case $\|DF\|_{\infty}=1$, $\gamma(n)$ is a linear function, so that 
${\tau}_*\geq \frac{1-e^{-1}}{C\|\nabla f\|_\infty}\ep^{-(\alpha\wedge 1)}$. 

In the coarse-grained version, if $\mu_F=1$ and $\|DF^n\|_\infty$ grows like in \eqref{power-law} 
with $\beta>0$,  the dissipation is slow:
$\tilde {\tau}_*\geq C\ep^{-(\alpha\wedge 1)/\beta}$.
In the case where $\|DF^n\|_\infty$ is uniformly bounded by some constant, 
the norm of the coarse-grained propagator
stays larger than some positive constant for all times, so that for small enough noise
$\tilde {\tau}_*$ is infinite. \hfill$\blacksquare$


\section{An upper bound of the dissipation time for mixing maps}
\label{UB}

For any two functions $f,\,h \in L^2_0(\IT^d)$, the dynamical correlation function for
the map $F$ is defined as the following function of $n\in\IN$ 
(see e.g. \cite{B}):
\bean
C_{f,h}(n)=C^0_{f,h}(n)=m(fU_{F}^{n}h)=\la \bar f,U_F^n h\ra= \la P_F^n \bar f,h\rangle.
\eean
The same quantity may be defined for the noisy evolution:
\bean
C^{\ep}_{f,h}(n)=m(fT_{\ep}^{n}h).
\eean
We recall that a map $F$ is {\em mixing} iff for any $f,\ h\in L^2_0$, 
\bean
C_{f,h}(n)\rightarrow 0, \quad \text{as}\quad n\to \infty.
\eean
The correlation function can easily be measured in (numerical or real-life) experiments, so it is
often used to characterize the dynamics of a system.

To focus the attention, we will only be concerned with maps 
for which correlations decay in a precise way. We assume that
there exist H\"older exponents $s_*,s\in \IR_+$, $0\leq s_*\leq s$
together with some decreasing
function $\Gamma(n)=\Gamma_{s_*,s}(n)$ with $\Gamma(n)\nto 0$, 
such that for any observables $f\in C_0^{s_*}(\IT^d)$,
$h\in C_0^{s}(\IT^d)$ 
and for sufficiently
small $\eps\geq 0$ (sometimes only for $\ep=0$), 
\bea
\label{dc}
\forall n\in\IN,\quad|C^{\ep}_{f,h}(n)|\leq \|f\|_{C^{s_*}} \|h\|_{C^{s}}\,\Gamma(n).
\eea
In general, such a bound can be proven only if the map $F$ has regularity 
$C^{s+1}$. The reason why we do not necessarily take the same norm for the functions 
$f$ and $h$ will be clear below.

We will be mainly interested in the following two types of decay

\begin{itemize}
\item[i)] Power-law decay: there exists $C>0$, $\beta>0$ such that,  
\bea
\Gamma(n)=C n^{-\beta}.
\eea
This behaviour is characteristic of intermittent maps, e.g. maps possessing one or several 
neutral orbits \cite{Bal-decay}.
\item[ii)] Exponential decay: there exists $C>0$, $0<\sigma<1$ such that, 
\bea
\Gamma(n)=C \sigma^n. 
\eea
Such a behaviour was proved in the case of uniformly expanding or hyperbolic maps on the torus 
(see Section~\ref{examples}), as well as many more general cases \cite{Bal-decay}.
\end{itemize}

The central result of this section is a relationship between, on one side, the
decay of noisy (resp. noiseless) correlations and on the other side, the small-noise
behaviour of the noisy (resp. coarse-graining) dissipation time. The intuitive picture is similar
to the one linking the local expansion rate to the dissipation: namely, a fast decay of 
correlations is generally due to the transition of ``soft'' into ``hard'' fluctuations of the
observable through the evolution, which is itself induced by large expansion rates of the map.
Still, as opposed to what we obtained in last Section, the following theorem and its corollary
yields {\em upper bounds} for the dissipation time.

\begin{thm}
\label{ccln}
Let $F$ be a volume preserving map on $\IT^d$ with correlations
decaying as in Eq.~(\ref{dc}) for some indices $s$, $s_*$ and decreasing function 
$\Gamma(n)$, at least in the noiseless limit $\ep=0$.
Assume that the noise generating function $g$ is $([s]+1)$-differentiable, and that
all its derivatives of order $|\bmal|_1\leq  [s]+1$ satisfy
\[|D^{\bmal}g(\bx)|\lesssim
\frac{1}{|\bx|^M},\quad |\bx|\gg 1,
\] with a power $M>d$.

Then there exist constants $\tilde C>0$, $\ep_o>0$ such that the coarse-grained propagator satisfies
\bea
\forall\ep\leq\ep_o,\quad\forall n\geq 0,\qquad
\|\tilde{T}^{(n)}_{\ep}\| \leq \tilde{C}\; \frac{\Gamma(n)}{\ep^{d+s+s_*}}.
\eea
If the decay of correlations \eqref{dc} also holds for sufficiently small $\ep>0$ 
(and assuming the Perron-Frobenius operator $P_F$ is bounded in $C^{s}(\IT^d)$), then
the noisy operator satisfies (for some constants $C>0$, $\ep_o>0$):
\bea
\forall \ep\leq\ep_o,\quad\forall n\geq 0,\qquad
\|T^{n}_{\ep}\| \leq C\; \frac{\Gamma(n)}{\ep^{d+s+s_*}}.
\eea
\end{thm}

From these estimates, we straightforwardly obtain the following bounds on both
dissipation times (the assumptions on $F$ and the noise generating function $g$
are the same as in the Theorem):

\begin{cor}
\label{correl-upper}

I) If the correlation function satisfies the bound \eqref{dc} for $\ep=0$, then
 the coarse-grained dissipation time is well defined
($\tilde{\tau}_*<\infty$). Moreover,
\begin{itemize}
\item[i)] if $\Gamma(n)\sim n^{-\beta}$ then there exists a constant
$\tilde{C}>0$ such that
\bean
\tilde{\tau}_* \leq \tilde{C} \ep^{-\frac{d+s+s_*}{\beta}}
\eean
\item[ii)] if  $\Gamma(n)\sim\sigma^n$ then there exists a constant
$\tilde{c}$ such that 
\bean
\tilde{\tau}_* \leq \frac{d+s+s_*}{|\ln\sigma|} \ln(\ep^{-1})+\tilde{c}, 
\eean
\end{itemize}

II) If Eq.~\eqref{dc} holds for sufficiently small $\ep>0$, then
\begin{itemize}
\item[i)] if $\Gamma(n)\sim n^{-\beta}$, there exists a constant $C>0$ such that
\bean
{\tau}_* \leq C \ep^{-\frac{d+s+s_*}{\beta}}
\eean
\item[ii)] if  $\Gamma(n)\sim\sigma^n$, there exists a constant $c$ such that
\bean
{\tau}_* \leq \frac{d+s+s_*}{|\ln\sigma|} \ln(\ep^{-1})+c. 
\eean
\end{itemize}
\end{cor}

\textbf{Proof of the Theorem.}

\underline{$1^{st}$ step}: We represent the action of $T_{\ep}^n$  
(resp. $\tilde{T}_{\ep}^{(n)}$) on an observable 
$f\in L^2_0(\IT^d)$
in terms of the correlation functions $C^\ep(n)$ (resp. $C(n)$). 
To do this we Fourier decompose both $T_{\ep}^{n+2}f$ and $f_1=U_F f$, and use Eq.~\eqref{specdecompo}:
\bean
T_{\ep}^{n+2}f &=& 
\sum_{0\neq \bj\in \IZ^d}\la \be_{\bj}, G_{\ep}U_{F}T_{\ep}^{n}G_{\ep}f_{1} \ra \be_{\bj}\\
&=&
\sum_{0\neq \bj\in \IZ^{d}}\sum_{0\neq \bk\in \IZ^{n}}\hat{f}_{1}(\bk)
\la G_{\ep}\be_{\bj},U_{F}T_{\ep}^{n}G_{\ep}\be_{\bk} \ra \be_{\bj}\\&=&
\sum_{0\neq 
\bj\in \IZ^{d}}\sum_{0\neq \bk\in \IZ^{d}}\hat{f}_{1}(\bk)
\hat{g}_{\ep}(\bj)  \hat{g}_{\ep}(\bk)
\la P_{F}\be_{\bj},T_{\ep}^{n}\be_{\bk} \ra 
  \be_{\bj}.
\eean
(remember that $\hat g$ is a real function). 
A similar computation for the coarse-grained propagator yields:
\bean
\tilde{T}_{\ep}^{(n)}f=\sum_{0\neq 
\bj\in \IZ^{d}}\sum_{0\neq \bk\in \IZ^{d}}\hat{f}(\bk)
\hat{g}_{\ep}(\bj)  \hat{g}_{\ep}(\bk)
\la \be_{\bj},U_{F}^{n}\be_{\bk} \ra 
  \be_{\bj}.
\eean
Taking the norms on both sides, we get in the noisy case:
\bea
\|T_{\ep}^{n+2}f\|^{2}&=&\sum_{0\neq \bj\in \IZ^{2}}
\bigg| \sum_{0\neq \bk\in \IZ^{d}}\hat{f}_{1}(\bk)
\la P_{F}\be_{\bj},T_{\ep}^{n}\be_{\bk} \ra \hat{g}_{\ep}(\bj) \hat{g}_{\ep}(\bk) \bigg|^{2}\nonumber\\
&\leq&
\sum_{0\neq \bj\in \IZ^{d}} 
\bigg( \sum_{0\neq \bk\in \IZ^{d}}|\hat{f}_{1}(\bk)|^{2}\bigg)
\sum_{0\neq \bk\in \IZ^{d}}|\la P_{F}\be_{\bj},T_{\ep}^{n}\be_{\bk} \ra|^{2}
|\hat{g}_{\ep}(\bj)\hat{g}_{\ep}(\bk)|^2\nonumber\\
\Longrightarrow \|T_{\ep}^{n+2}f\|^{2}&\leq& \|f_1\|^2
\sum_{0\not = \bj,\bk\in \IZ^{d}}|C^{\ep}_{P_{F}\be_{-\bj},\be_{\bk}}(n)|^{2}
  |\hat{g}(\ep\bj)\hat{g}(\ep\bk)|^2, \label{noisy-ineq} 
\eea
and in the coarse-graining case
\bea
\label{coarse-ineq}
\|\tilde{T}_{\ep}^{(n)}f\|^2 \leq \|f\|^2 \sum_{0\neq \bj,\bk\in \IZ^{d}}|C_{\be_{-\bj},\be_{\bk}}(n)|^{2}
  |\hat{g}(\ep\bj)\hat{g}(\ep\bk)|^2. 
\eea
These two expressions explicitly relate the dissipation with the correlation functions.
\medskip

\underline{$2^{nd}$ step:} 
We now apply the estimates (\ref{dc}) on correlations for the observables $\be_{\bk}$, $\be_{-\bj}$, 
$P_{F}\be_{-\bj}$. 
In the coarse-grained case, it yields (using simple bounds of the type of Eq.~\eqref{cosin}):
\bean
\forall\bj,\ \bk\in\IZ^d\setminus \{0\},\qquad |C_{\be_{-\bj},\be_{\bk}}(n)|
\leq C'\;|\bj|^{s}|\bk|^{s_*}\Gamma(n).
\eean
In the noisy case, we need to assume that the Perron-Frobenius operator $P_F$ is bounded in the space
$C^{s}(\IT^d)$. This property is in general a prerequisite in the proof of estimates
of the type \eqref{dc}, so this assumption is quite natural here.
\bea
\forall\bj,\ \bk\in\IZ^d\setminus \{0\},\qquad
|C^{\ep}_{P_{F}\be_{-\bj},\be_{\bk}}(n)| &\leq& C\; \|P_{F}\be_{-\bj}\|_{C^{s}} \|
 \be_{\bk}\|_{C^{s_*}}\Gamma(n)\non\\
&\leq& C\|P_F\|_{C^{s}} |\bj|^{s}|\bk|^{s_*}\Gamma(n).\label{Cep}
\eea
We insert these bounds on the decay of correlations in the expressions (\ref{noisy-ineq}-\ref{coarse-ineq}), 
for instance in the coarse-grained case we get:
\bea
\forall n\geq 0,\qquad
\|\tilde{T}^{(n)}_{\ep}\|^{2} \leq C\; \Gamma(n)^2 \bigg(\ep^{-(s+s_*)}\sum_{0\neq \bk\in \IZ^{d}}|\ep\bk|^{s+s_*}
\hat{g}(\ep\bk)^2 \bigg)^{2}.\label{step2}
\eea

\medskip

\underline{$3^{rd}$ step:}
We finally estimate the $\ep$-dependence of the RHS of the above inequality. 
Up to a factor $\ep^{-d}$, the sum in the brackets is
a Riemann sum for the integral $\int |\bxi|^{s+s_*}\hat g(\bxi)^2 d\bxi<\infty$. 
A precise connection is given in the following lemma, proven in Appendix~\ref{app4}:
\begin{lem}
\label{sum-int}
Let $f\in C^0(\IR^d)$ be symmetric w.r.t. the origin and decaying at infinity as $|f(\bx)|\lesssim |\bx|^{-M}$ with $M>d$. 
Then the following estimate holds in the limit $\ep\to 0$:
\bea
\label{sum-integ}
\ep^d\sum_{\bk\in\IZ^d}\hat f(\ep\bk)^2
=\int_{\IR^{d}} \hat f(\bxi)^2 d\bxi+\cO(\ep^M).
\eea
\end{lem}
Let $m\in\IN$ satisfy $2m\leq s+s_*\leq 2m+2$ (notice that $m\leq [s]$ since we assumed $s_*\leq s$).
From the obvious inequality 
$$
\forall x>0,\qquad x^{s+s_*}\leq x^{2m} +x^{2m+2},
$$ 
we may replace in the RHS of \eqref{step2} the factor $|\ep\bk|^{s+s_*}$ by $|\ep\bk|^{2m}+|\ep\bk|^{2m+2}$. 
Applying Lemma~\ref{sum-int} to the derivatives of $g$ of order $m$ and $m+1$, we end up with
the following upper bound, which proves the first part of the
theorem: 
\bean
\|\tilde{T}^{(n)}_{\ep}\|^{2} 
&\leq&C\;\Gamma(n)^2 
\left(\frac{1}{\ep^{d+s+s_*}}\int_{\IR^{d}} \big(|\bxi|^{2m}+|\bxi|^{2(m+1)}\big)\hat g(\bxi)^2 d\bxi+\cO(\ep^M)\right)^{2} \\
&\leq& C'\;\frac{\Gamma(n)^2}{\ep^{2(d+s+s_*)}}\;\|g\|^4_{H^{m+1}}.
\eean
The computations follow identically
for the case of the noisy operator, yielding the second part of the theorem. \hfill$\blacksquare$


\section{Some examples of mixing maps}\label{examples}

In this section we apply the results of the last two sections to several classes of 
Lebesgue measure preserving maps which have been proven to be mixing, with various types of correlation
decays. 

\medskip 

\textbf{Remark:}
In general, an expanding or hyperbolic map on the torus does not preserve the Lebesgue measure, so the first
step is to precise with respect to which invariant measure 
one wants to study the ergodic properties. In the ``nice'' cases, one can prove the existence 
and uniqueness of a ``physical measure'', which is ergodic for the map $F$, and then study the (noisy or
noiseless) mixing properties with respect to this measure. As pointed out in the introduction, in the present paper 
we only consider maps for which the physical measure is the Lebesgue measure.


\subsection{Decays of correlations}
\label{noiseless}
We briefly summarize the results concerning mixing maps, mostly relying on the
review \cite{Bal-decay} and the book \cite{B}.

Let us first consider the noiseless correlations. A common
route to prove that a map $F$ is {\sl exponentially} mixing consists in identifying an
invariant space $\cB$ of densities, 
on which the Perron-Frobenius has
the following spectral structure: except for the eigenvalue $1$ associated with the constant
density, the rest of the spectrum (on the subspace $\cB_0$ of zero-mean densities)
is contained inside a disk of radius $\sigma<1$ centered
at the origin. More precisely, the spectrum on $\cB_0$ may consist in isolated eigenvalues (called {\sl resonances})
$\{\lambda_i\}$ and of some essential spectrum in the disk of radius $r_{ess}<\sigma$.
$P_F$ is said to be {\sl quasicompact} on $\cB$.
 This spectral structure implies that there is some constant $C>0$ such that
for any $f\in\cB_0$, $h\in\cB^*_0$,
\bea
\label{expon1}
\forall n\geq 0,\qquad |C_{f,h}(n)|=|\la h,P_F^n f\ra_{\cB^*,\cB}|
\leq \|h\|_{\cB^*}\|P_F^n\|_{\cB_0}\|f\|_{\cB}
\leq C\,\|h\|_{\cB^*} \|f\|_{\cB}\,\sigma^n.
\eea

For the maps we study, the spectrum of $P_F$ on $L^2_0$ 
intersects the unit circle, so $\cB$ cannot simply be the Hilbert space $L^2$.
Depending on the properties of the map, $\cB$ can be a Fr\'echet space of analytic functions,
a Banach space of bounded variation, H\"older or $C^{s}$ functions 
(see Section~\ref{expanding}); it may also
be a space of generalized functions lying outside $L^2$ (see Section~\ref{Anosov}). 
No matter how complicated $\cB$ can be, in general there exist H\"older exponents 
$0\leq s_*\leq s$ such that $C^s$  (resp. $C^{s_*}$) 
embeds continuously in $\cB$ (resp. in its dual $\cB^*$). As a result, the upper bound \eqref{expon1} can
be specialized to functions $f\in C_0^s$, $h\in C_0^{s_*}$ as follows: 
\bea
\label{expon}
\forall n\geq 0,\qquad |C_{f,h}(n)|\leq C\,\|h\|_{C^{s_*}} \|f\|_{C^{s}}\,\sigma^n.
\eea
This is the form of upper bound we used in Theorem~\ref{ccln}.

This strategy of proof has been applied to several types of maps, including
the (noninvertible) expanding maps and the Anosov or Axiom-A diffeomorphisms on a
compact manifold.  Exponential decay of correlations has also been proven (using various methods) 
for piecewise expanding maps on the
interval, some nonuniformly hyperbolic/expanding maps,  
some expanding or hyperbolic maps with singularities. 

Other types of decay occur as well: for instance, 
a polynomial decay of correlations $C_{f,h}(n)\lesssim n^{-\beta}$ was
shown to be optimal for some ``intermittent'' systems, like a
one-dimensional map expanding everywhere except at a fixed ``neutral'' point (such maps are sometimes
called ``almost expanding'' or ``almost hyperbolic'').

\medskip

There exist fewer results on the decay of correlations for stochastic perturbations of deterministic maps, 
like our noisy evolution $T_\ep$. In general, one wants to prove
\emph{strong stochastic stability}, that is stability of the invariant measure and of 
the rate of decay of the correlations in the small-noise limit. In our case, only the
second point needs to be proven, since the Lebesgue measure remains invariant after switching
on the noise.

Stochastic stability has been proven for smooth uniformly expanding
maps \cite{BalYou} (see next subsection), some nonuniformly expanding or piecewise expanding maps.
It has been shown also for uniformly hyperbolic (Anosov) maps on the $2$-dimensional torus \cite{BKL}
(see subection~\ref{Anosov}). In all those cases, the mixing is exponential, 
so the stability of the decay \eqref{expon} means that
for small enough $\ep>0$, there exists a radius $\sigma_\ep\epto \sigma$ such that for any 
$f\in C^{s}_0$, $h\in C^{s_*}_0$,
\bea
\label{noisy-expon}
\forall n>0,\qquad |C^\ep_{f,h}(n)|\leq C\,\|h\|_{C^{s_*}}\,\|f\|_{C^{s}}\,\sigma_\ep^n.
\eea
Here the constant $C>0$ can be taken independent of $\ep$.

In the next two sections, we describe in more detail the cases
of smooth uniformly expanding maps and Anosov diffeomorphisms on the torus.


\subsection{Smooth uniformly expanding maps}
\label{expanding}
Let $F$ be a $C^{s+1}$ map on $\IT^d$ (with $s\geq 0$). 
Assume there exists $\lam>1$ such that for any $\bx\in\IT^d$ and any $\bv$ in the
tangent space $T_{\bx}\IT^d$, $\|DF(\bx)\bv\|\geq\lam \|\bv\|$ (we assume that $\lam$ is the
largest such constant). Such a map is called
uniformly expanding. In general, it admits a unique absolutely continuous invariant
probability measure; here we restrict ourselves to maps for which this measure is the Lebesgue measure.

Ruelle \cite{Ru89} proved that the Perron-Frobenius operator $P_F$ of such a map is quasicompact on the space
$C^s(\IT^d)$, and that its essential spectrum is contained inside the disk of radius $r_1=\lam^{-s}$. In general, 
one has little information on the possible discrete spectrum outside this disk (upper bounds on the
decay rate have been obtained in the case of an expanding map of 
regularity $C^{1+\eta}$ \cite{B}).  
Strong stochastic stability for such maps was proven in \cite{BalYou}, with a 
more general definition of the noise than the one we gave.

For all these cases, one can take $s_*=0$, since the 
continuous functions are continuously embedded in any space $(C^s)^*$. 

\subsection*{Case of a linear expanding map}
We describe the simplest example possible for such a map, namely the angle-doubling map on $\IT^1$ defined as
$F(x)=2\,x\bmod 1$. This map is real analytic, with uniform expansion rate $\lam=2$. Due to its linearity, the
dynamics of this map (as well as its noisy version) 
is simple to express in the basis of Fourier modes $\be_k(x)=e^{2i\pi kx}$:
\bean
\forall k\in\IZ,\quad U_{F}\be_{k}&=&\be_{2k}\\
\Longrightarrow T_\ep\be_{k}&=&\hat g(\ep k)\be_{2k}\\
\Longrightarrow T^n_\ep\be_{k}&=&\big[\prod_{j=1}^{n}\hat g(\ep\, 2^j k)\big]\be_{2^n k}
\eean
The computation is even simpler for the coarse-grained propagator:
$$
\tilde T^{(n)}_\ep \be_{k}=\hat g(k)\hat g(2^n k)\be_{2^n k}.
$$
To fix the ideas, we consider the $\alpha$-stable noise $\hat g(\xi)=e^{-|\xi|^\alpha}$ for some $0<\alpha\leq 2$.
One easily checks that for any $n\geq 1$, 
\bean
\|T_\ep^n\|&=&\|T_\ep^n \be_1\|=\exp\bigg\{-\ep^\alpha \frac{2^{n\alpha}-1}{1-2^{-\alpha}}\bigg\},\\
\|\tilde{T}^{(n)}_\ep\|&=&\exp\{-\ep^\alpha (2^{n\alpha}+1)\}.
\eean
For any $\ep>0$, these decays are super-exponential: the spectrum of $T_\ep$ on $L^2_0$
is reduced to $\{0\}$ for any $\ep>0$ (the spectrum of $U_F$ is the full unit disk).
From this explicit expression, we get both dissipation times:
\bea
\label{exact-asympt}
{\tau}_* =\frac{1}{\ln 2}\ln(\ep^{-1})+\cO(1),
\qquad \tilde{\tau}_*=\frac{1}{\ln 2}\ln(\ep^{-1})+\cO(1).
\eea
For this linear map, $\|DF\|_\infty=\mu_F=2$, so this estimate is in agreement
with the lower bounds \eqref{nlb}, the latter being sharp if $\alpha\in [1,2]$.
On the other hand, $\ln 2$ is also equal with the Kolmogorov-Sinai (K-S) entropy $h(F)$ of $F$.
Therefore, for this linear map the dissipation rate constant exists, and is equal to $\frac{1}{h(F)}$.

To compare these exact asymptotics with the upper bounds of Corollary~\ref{correl-upper}, we estimate the
correlation functions $C_{f,h}(n)$ on the spaces $C^{s}(\IT^1)$. We give below a short proof in the
case $s>\half$. We will use the following Fourier estimates \cite{Zyg}:
$$
\exists C>0,\quad,\forall f\in C_0^{s}(\IT^1),\quad \forall k\neq 0,\qquad 
|\hat f(k)|\leq C\frac{\|f\|_{C^{s}}}{|k|^{s}}.
$$
Therefore, writing the correlation function as a Fourier series, we get:
\bea
\|P_F^n f\|^2&=&\sum_{0\not =k\in \IZ} |\hat f(2^n k)|^2
\leq\sum_{0\not = k\in \IZ}\lp C\frac{\|f\|_{C^{s}}}{|2^n k|^{s}}\rp^2\non\\
\Longrightarrow \|P_F^n f\|&\leq&C'\,\frac{\|f\|_{C^{s}}}{(2^{s})^n}.\label{angle-d-spectrum}
\eea
This estimate yields a decay of the correlation function 
as in Eq.~\eqref{expon}, with a rate $\sigma=2^{-s}$ and $s_*=0$. One can check that this rate is
sharp for functions in $C^{s}$: indeed, any $z\in\IC$, $|z|<2^{-s}$ is an eigenvalue
of $P_F$ on that space.
Applying the Corollary~\ref{correl-upper},{\em I ii)}, we
get that for any $s>1/2$, there exists a constant $\tilde{c}$ such that
\bea
\label{angle-d-upper}
\tilde{\tau}_*\leq \frac{1+s}{s\ln 2}\ln(\ep^{-1})+\tilde{c}
\eea
for sufficiently small $\ep$. The exact dissipation rate
constant $1/\ln 2$ is recovered only for large $s$.

A straightforward computation shows that the estimate
\eqref{angle-d-spectrum} also holds if one replaces $P_F$ by $P_F\circ G_\ep$; hence the noisy
correlation function dynamics satisfies the same uniform  upper bound as the noiseless one, with the
decay rate
$\sigma_\ep=2^{-s}$. As a result, the upper bound on ${\tau}_*$ given by 
Corollary~\ref{correl-upper},{\em II ii)} is the same as in  Eq.~\eqref{angle-d-upper}.


\subsection{Anosov diffeomorphisms on the torus}
\label{Anosov}
We recall that a diffeomorphism $F:\IT^d\mapsto \IT^d$ is called
Anosov if it is uniformly hyperbolic:
there exist constants $A>0$ and $0<\lam_{s}<1<\lam_{u}$ such that
at each $\bx\in\IT^d$ the tangent space $T_{\bx}\IT^d$
admits the direct sum decomposition $T_{\bx}\IT^d= E^{s}_{\bx} \oplus E^{u}_{\bx}$ 
into stable and unstable subspaces such that for every $n\in \IN$,
\bean
(D_{\bx}F)(E^{s}_{\bx})=E^{s}_{F\bx}, \qquad \|(D_{\bx}F^n)_{|E^{s}_{\bx}}\|\leq A \lam_{s}^n;\\
(D_{\bx}F)(E^{u}_{\bx})=E^{u}_{F\bx}, \qquad \|(D_{\bx}F^{-n})_{|E^{u}_{\bx}}\|\leq A \lam_{u}^{-n}.
\eean
These inequalities have obvious consequences on the expansion rates of $F$ and $F^{-1}$, for 
instance they imply $\|DF\|_\infty^n\geq\|DF^n\|_\infty\geq A^{-1}\lam_u^n$. As a consequence,
the quantities of interest in Theorem~\ref{nln}, {\em i)} satisfy
\bean
\|DF\|_\infty\wedge\|DF^{-1}\|_\infty&\geq& \lam_u\wedge\lam_s^{-1},\\
\mu_F\wedge\mu_{F^{-1}}&\geq& \lam_u\wedge\lam_s^{-1}.
\eean 
All these expansion rates are $>1$, so both noisy and coarse-grained dissipation times admit
logarithmic lower bounds as in Eq.~\eqref{nlb}.

Exponential mixing has been proven for Anosov diffeomorphisms of regularity $C^{1+\eta}$ ($0<\eta<1$)
by Bowen \cite{Bow}, using
symbolic dynamics; the exponential decay is then valid for H\"older 
observables in $C^{\eta'}$ for some $0<\eta'<\eta$.
Because we are also interested in the noisy dynamics, we will
refer to a more recent work \cite{BKL} concerning $C^3$ Anosov maps on $\IT^d$, which bypasses
symbolic dynamics. The authors
construct an invariant Banach space $\cB$ of generalized functions on the phase space, 
such that the Perron-Frobenius operator is quasicompact on this space. One subtlety (compared with the
case of expanding maps) is that $\cB$ explicitly depends
on the (un)stable foliations of the map $F$ on $\IT^d$. Vaguely speaking, the elements of $\cB$ are
``smooth'' along this unstable direction $E^{u}_{\bx}$, but can be singular (``dual of smooth'') 
along the stable foliation $E^{s}_{\bx}$. 
The space $\cB$ is
the completion of $C^1(\IT^d)$ with respect to a norm $\|\cdot\|_{\cB}$ adapted to these
foliations. 
This norm, and therefore $\cB$, are defined in terms of a parameter $0<\beta<1$, the choice of which
depends on the {\sl regularity} of the unstable foliation. 
In general, the latter is 
$\tau$-H\"older continuous, for some exponent $0<\tau<2$. Then, the authors prove that if
one takes $\beta<\tau\wedge 1$, then 
the essential spectrum of $P_F$ on $\cB$ has a radius smaller than 
$r_\beta=\max(\lambda_u^{-1},\lambda_s^{\beta})$. This upper bound
is sharper if $\beta$ can be taken close to $1$, that is, if the foliation is $C^1$. 
This is the case for smooth area-preserving Anosov maps on $\IT^2$, for which the foliations have regularity
$C^{2-\delta}$ for any $\delta>0$ \cite{HK}. 
The operator $P_F$ may have isolated eigenvalues (resonances) $1>|\lambda_i|>r_\beta$, 
corresponding to eigenstates in $\cB_0$ 
which are genuine distributions $\not\in L^2_0$. There is (to our knowledge) no simple
general upper bound for the largest resonance $|\lambda_1|$ 
in terms of the expansion parameters $(\lambda_u,\ \lambda_s)$. 

By construction, the space $C^1(\IT^d)$ embeds continuously in both $\cB$ and its dual $\cB^*$, 
so that one can take $s=s_*=1$ in Eq.~\eqref{expon}. Therefore, for any
$\sigma_\beta>\max(|\lam_1|,r_\beta)$, there is some constant $C>0$ such that for any $f$, $h\in C_0^1(\IT^d)$,
\bea
\label{corr-anosov}
\forall n>0,\qquad |C_{f,h}(n)|\leq C\,\|h\|_{C^1}\,\|f\|_{C^1}\,\sigma_\beta^n.
\eea
In the proof of Theorem~\ref{ccln} (Step 3), for the case $s=s_*=1$ we only need to assume that
the noise generating function $g$ is $C^1$ with fast-decaying first derivatives. The fast decay
implies that the first moment of $g$ is finite (that is, one can take $\alpha\geq 1$).

The noisy propagator $G_\ep P_{F}$ is also analyzed in \cite{BKL}.
If the unstable foliation has regularity $C^{1+\eta}$ with
$\eta>0$ (for instance for any $C^3$  Anosov 
diffeomorphism on $\IT^{2}$), and assuming that the noise generating
function $g\in C^2(\IR^d)$ has {\em compact support}\footnote{The condition of compact support 
could probably be relaxed to one of fast decrease
at infinity (C. Liverani, private communication)}, the authors prove
the strong spectral stability of the Perron-Frobenius operator $P_F$ on any space $\cB$ defined
with a parameter $\beta'<\eta$. Notice that this constraint on $\beta'$ is stronger than
in the noiseless case, where one could take any $\beta<1$; the spectral radius $\sigma_{\beta'}$
may accordingly be larger than $\sigma_\beta$. 
Modulo the replacement of $\sigma_\beta$ by $\sigma_{\beta'}$, the estimate \eqref{corr-anosov} 
therefore applies to
the noisy correlation function $C^{\ep}_{f,h}(n)$ as long as $\ep$ is small enough.

We collect below the results regarding the dissipation 
time of $C^3$ Anosov maps on the torus.

\begin{thm}
\label{TAnosov}
Let $F$ be a volume preserving $C^{3}$ Anosov diffeomorphisms on $\IT^d$,
and let the noise generating function be $C^1$ with fast decay at infinity.

I) Then there exist 
$\mu\geq\lam_u\wedge\lam_s^{-1}$, $0<\tilde{\sigma}<1$ and
 $\tilde{C}>0$ such that
the dissipation time of the coarse-grained dynamics satisfies
\bean
  \frac{1}{\ln \mu }\ln(\ep^{-1}) -\tilde{C} \leq
  \tilde{\tau}_* \leq \frac{d+2}{|\ln \tilde{\sigma}|} \ln(\ep^{-1}) + \tilde{C},
\eean

II) If in addition $F$ has $C^{1+\eta}$-regular foliations, 
and $g\in C^2(\IR^{d})$ is compactly 
supported, then there exist $\tilde{\sigma}\leq\sigma<1$ and  $C$ such that the 
dissipation time of the noisy dynamics satisfies
\bean
  \frac{1}{\ln \|DF\|_\infty}\ln(\ep^{-1}) - C
\leq  {\tau}_* \leq \frac{d+2}{|\ln \sigma|}\ln(\ep^{-1}) + C
\eean
\end{thm}


\subsection{Ergodic linear automorphisms of the torus}
\label{ergodicaut}
In this section we describe examples of Anosov maps for which the dissipation
time can be precisely determined. One can even compute
the dissipation rate constant; we will mention the
connection of the latter with the Kolmogorov-Sinai entropy. 
After the simple example of
a (generalized) cat map on the $2$-dimensional torus \cite{Arnold}, we recall 
the results obtained in \cite{FW} for $d$-dimensional ergodic automorphisms.

Throughout this section, we take a noise generating function of the type \eqref{Akernel}
for a certain $\alpha\in (0,2]$, with $\bQ=I$:
\bea
\label{ak}
g_{\ep,\alpha}(\bx):=\sum_{\bk \in \IZ^{d}}e^{-|\ep \bk|^{\alpha}}\be_{\bk}(\bx).
\eea
Notice that this function is not compactly supported, therefore it does not satisfy {\sl stricto sensu}
the assumptions required in \cite{BKL} to prove the strong spectral stability of $P_F$ (see the footnote
after Theorem~\ref{TAnosov}).


\subsubsection*{$2$-dimensional cat map}
We take $A\in SL(2,\IZ)$ 
with $|TrA|>2$ and define the dynamics as $F(\bx)=F_A(\bx)=A^t\bx  \mod 1$ ($A^t$ is the transposed
matrix of $A$). This map
is of Anosov type, and preserves the Lebesgue measure.
The dynamics is easy to express on the Fourier modes:
\bean
U_{F}\,\be_{\bk}&=&\be_{A\bk}\\
T^n_{\ep}\,\be_{\bk}&=&e^{-\sum_{l=1}^n|\ep A^l \bk|^{\alpha}}\be_{A^n \bk} \\
\Longrightarrow \|T^n_{\ep}\|&=& \exp\lp-\ep^{\alpha}\min_{0\not=\bk\in\IZ^2} 
\sum_{l=1}^n|A^l \bk|^{\alpha}\rp.
\eean
Similarly, 
$$
\|\tilde{T}^{(n)}_{\ep}\|=\exp\lp-\ep^{\alpha}
\min_{0\not=\bk\in\IZ^2}(|\bk|^{\alpha}+ |A^n \bk|^{\alpha})\rp.
$$
Let us call $\lam$, $\lam^{-1}$ the eigenvalues of $A$, with the convention $|\lam|>1$. 
One can easily show that there are two constants $0<C_1<C_2$ such that for any $n>0$,
\bean
C_1\,|\lam|^{n\alpha/2}&\leq& \min_{0\neq\bk\in \IZ^2} \sum_{l=1}^{n}|A^{l}\bk|^{\alpha}
\leq C_2\,|\lam|^{n\alpha/2},\\
C_1\,|\lam|^{n\alpha/2}&\leq& \min_{0\neq\bk\in \IZ^2} (|\bk|^{\alpha}+|A^{n}\bk|^{\alpha})
\leq C_2\,|\lam|^{n\alpha/2}.
\eean
The above estimate yields the following asymptotics in the small-$\ep$ limit:
\bean
 {\tau}_*= \frac{2}{\ln|\lam|}\ln(\ep^{-1})+\cO(1),
\qquad\tilde{\tau}_*= \frac{2}{\ln|\lam|}\ln(\ep^{-1})+\cO(1).
\eean
As in the case of linear expanding maps, the dissipation rate constant exists, and 
seems related with the K-S entropy of the linear map $h(F)=\ln|\lam|$.

Let us compare these exact asymptotics with the bounds obtained in the previous sections.
The cat map being linear, one has 
$\|DF\|_\infty=\|A^t\|\geq |\lam|$,  $\|DF^{-1}\|_\infty=\|(A^t)^{-1}\|\geq |\lam|$ 
(with equality iff $A$ is symmetric). Since $A$ is diagonalizable, we have 
$\|DF^n\|_\infty=\|(A^t)^n\|\sim |\lam|^n$, 
so that $\mu_F=\mu_{F^{-1}}=|\lam|$. Therefore, the lower bounds for the dissipation times given in 
Theorem~\ref{nln}{\em i)} are strictly smaller than the exact rates derived above, they differ from
the latter by a factor $\leq 1/2$. 

To estimate the dynamical correlations, we use Fourier decomposition and proceed as in \cite{CrawCar}:
we construct a ``primitive subset'' $\cS$ of the Fourier plane $\IZ^2\setminus \{0\}$,
such that each Fourier orbit intersects $\cS$ once and only once. This subset looks like the union of
$4$ angular sectors of the plane, and has the following crucial property:
\bea
\label{est1}
\exists c>0,\quad\forall \bk\in\cS,\quad\forall n\in \IZ,\quad |A^n\bk|\geq c|\lam|^{|n|}|\bk|.
\eea
The correlation function between two observables $f,\ h\in L^2_0(\IT^d)$ may be written as: 
$$
C_{f,h}(n)=\sum_{0\neq \bk\in\IZ^2}\hat f(-\bk)\hat h(A^n \bk)=
\sum_{\bk\in\cS}\sum_{l\in\IZ}\hat f(-A^l \bk)\hat h(A^{l+n} \bk).
$$
Let us assume $s>1$, and take $f\in C^{s}_0(\IT^2)$: its Fourier coefficients decrease 
for all $0\neq\bk\in\IZ^2$ as
$$
 |\hat f(\bk)|\leq C \frac{\|f\|_{C^{s}}}{|\bk|^{s}}.
$$
Using this decrease as well as the estimate \eqref{est1}, we find that 
$$
|C_{f,h}(n)|\leq C^2\,\|f\|_{C^{s}}\,\|h\|_{C^{s}}
\sum_{\bk\in\cS}\sum_{l\in\IZ} \frac{1}{(c^2|\lam|^{|l+n|+|l|}|k|^2)^{s}}.
$$
The sum over $\bk$ converges, and the
sum over $l$ is bounded from above by $\frac{C}{|\lam|^{ns}}$. 
This implies that for any $s>1$, the correlation function decreases as 
in Eq.~\eqref{expon}, with
the rate  $\sigma=|\lambda|^{-s}$. Actually,
this decrease holds as well in the case $0<s\leq 1$,
but the proof is different. 
From there, Corollary~\ref{correl-upper}{\em I ii)} yields the upper bound
\bea
\tilde{\tau}_*\leq \frac{2+2s}{s\ln|\lam|}\ln(\ep^{-1})+c.
\eea

This method can be straightforwardly adapted to prove that the noisy correlation function 
$C^\ep_{f,h}(n)$ decreases as fast as the
noiseless one (in case $s>1$). 
As a result, we obtain the same upper bound for ${\tau}_*$ as for $\tilde{\tau}_*$. 
We notice that, as in the case of the
angle-doubling map, the constant in the upper bound converges to the exact rate only
in the limit $s\gg 1$, while for $s=1$ 
(cf. the discussion preceding Thm.~\ref{TAnosov}), the constant is twice larger.


\subsubsection*{Toral automorphisms in higher dimensions}
We consider toral automorphisms $F=F_A$ given by matrices $A\in SL(d,\IZ)$; 
such an automorphism is ergodic
iff none of the eigenvalues of $A$ are roots of unity.
The K-S entropy $h(F)$ of $F$ is given by the formula \cite{Yuz}
\bea
\label{entropy}
 h(F)=\sum_{|\lam_{j}| > 1}\ln{|\lam_{j}|},
\eea 
where $\lam_{j}$ are the eigenvalues of $F$. All ergodic toral automorphisms have positive entropy.
We denote by $P$ the characteristic polynomial of the matrix $A$ and by 
$\{P_{1},...,P_{r}\}$ the complete set of its distinct irreducible factors (over $\IQ$). 
Let $d_{j}$ denote the degree of the polynomial $P_{j}$ and $h_{j}$ the   
K-S entropy of a toral automorphism with the characteristic polynomial $P_{j}$. 
For each $P_{j}$ we define its dimensionally averaged K-S entropy as 
\bea
 \hat{h}_{j}=\frac{h_{j}}{d_{j}}.
\eea 
For the whole matrix $F$ we define its minimal dimensionally averaged entropy (denoted $\hat{h}(F)$) as
\bean
 \hat{h}(F)=\min_{j=1,...,r}\hat{h}_{j} 
\eean
Notice that for an {\em irreducible} matrix $A$, that is a matrix admitting no proper invariant
rational subspace, this quantity reduces to $\frac{h(F)}{d}$.

Our notation regarding the noise parameter 
slightly differs from the one used in \cite{FW} (one has to replace 
$\vep^{2\alpha}$ by $\ep$); the results obtained there
read as follows:

\begin{prop}\cite{FW}
\label{p}
Let $F=F_A$ be a toral automorphism on $\IT^d$, and assume the noise kernel to be
as in Eq.~(\ref{ak}). Then, in the limit of small noise,
 
i) both dissipation times have a power-law behaviour iff $F$ is not ergodic.

ii) both dissipation times have logarithmic behaviour iff $F$ is ergodic. 

iii) if $F$ is ergodic and $A$ is diagonalizable then 
\bean
{\tau}_*\approx\tilde{\tau}_* \approx \frac{1}{\hat{h}(F)}\ln(\ep^{-1}).
\eean
\end{prop}

We end this section with a remark that the small-noise asymptotics of the dissipation time is
insensitive to a super-exponential decay of the correlation functions. We illustrate
this fact by the following result on the decay of correlations 
for $d$-dimensional toral automorphisms, which can be proven along the same lines
as the above proposition:

\begin{prop}
\label{supexp}
Let $F$ be a diagonalizable ergodic toral automorphism and $\lam$ any constant such that 
 $0<\lam<\hat{h}(F)$. 
Then for any $f,h \in L^{2}_{0}(\IT^{2d})$
\bean
C^{\ep}_{f,h}(n)\leq \|f\|\,\|h\|\,e^{-\ep^2 \lam^n}
\eean
Let $f,\,h \in G_{\ep}(L^{2}_{0}(\IT^{2d}))$ be smooth observables. Then
\bean
C_{f,h}(n)\leq \|G_{\ep}^{-1}f\|\,\|G_{\ep}^{-1}h\|\,e^{-\ep^2 \lam^n}.
\eean
\end{prop}


\section{Conclusion}
We have investigated the effect of noise or coarse-graining on the
dynamics generated by a conservative map, in particular the connection between 
the speed of dissipation of the noisy dynamics and the
spectral and dynamical properties of the underlying map.

We restricted ourselves to the case of volume preserving
maps on the $d$-dimensional torus. The choice of the torus allowed us
to use the Fourier transformation in our proofs, that is, harmonic analysis 
on that manifold. Most of our results can certainly be generalized to volume preserving
maps on other compact Riemannian manifolds. The choice of the torus was also guided by
the existence of simple volume preserving Anosov maps on it, most notably the linear
examples presented in Sections~\ref{expanding}, \ref{ergodicaut}. As explained at the beginning of Section~\ref{examples},
one may also want to extend the results to maps which do not leave invariant the Lebesgue measure, 
but still admit a `physical measure', as is the case for uniformly expanding or
Anosov maps (the physical measure is of SRB type, that is, its projection along
the unstable manifold is absolutely continuous \cite{B}). The noisy 
perturbations of these maps have been analyzed as well \cite{BalYou,BKL}, in particular
the strong spectral stability of the Perron-Frobenius operator holds under
the same conditions as for the volume preserving maps. Although the equilibrium measure
is more complicated, it might be possible to define and study a dissipation time 
in this more general framework.


As we explained in Corollary~\ref{weakmix}, the dissipation of a
non-weakly-mixing map is governed only by the nontrivial
eigenstates of the Koopman operator, the speed of dissipation depending on the smoothness
of these eigenstates. The asymptotics of the noisy dissipation time is generally
a power law in $\ep^{-1}$, while there
remains a possibility of faster dissipation if all eigenstates are singular enough
(see Remarks after Corollary \ref{weakmix}).

The more interesting results concern the mixing dynamics, in particular when the mixing occurs
exponentially fast.
We then proved that both noisy and coarse-graining dissipation times
behaved as the logarithm of $\eps^{-1}$ in the small noise limit. This dependence can be understood
through the time evolution in Fourier space: a mixing map transforms long wavelength fluctuations
into short wavelength ones, the latter being damped fast by the noise operator. This evolution in Fourier space
is typical of uniformly expanding/hyperbolic maps; it is already responsible for the fast decay of 
dynamical correlations. The link between this decay and the dissipation is explicited in Theorem~\ref{ccln},
and its application to Anosov maps is given in Theorem~\ref{TAnosov}.

In this context, we were unable to solve the problem of the
existence and value of the dissipation rate constant (that is, the prefactor in
front of $\ln(\ep^{-1})$). We obtained lower, resp. upper bounds for this
constant, in terms of the local expanding rate, resp. the rate of decay of correlations. 
It is not clear whether this constant can in general be related with the measure-theoretic (KS) 
entropy of the map, as is the case for linear automorphisms (modulo some algebraic subtleties \cite{FW}).

Although we used the spectral estimates of Theorem~\ref{gb} mostly for the case of non-weakly-mixing
maps, it also makes sense to use that theorem in the reverse direction in the case of mixing maps, that
is, deduce (pseudo)spectral properties of the noisy propagators, starting from the dissipation time estimates
obtained in Corollary~\ref{correl-upper}. The analysis of the pseudospectrum of $T_\ep$ for mixing maps
would complement the spectral one \cite{BKL,N}. We did not enter this aspect in the main text, 
because our attention was devoted to obtaining information on the dissipation time.

The results of the present paper concern
classical (i.e. non-quantum) dynamical systems.
The quantization of both linear and nonlinear maps on a symplectic
(even-dimensional) torus as a phase space has been studied in a number of 
works \cite{HB,BNS,tabor,KMR,GraDE,Z}. Several recent studies deal with some form of 
noise, or decoherence, in discrete-time quantum dynamics \cite{braun,fishman,manderfeld,saraceno,N}.
While most of these works concentrate on spectral or entropic
properties of noiseless/noisy dynamics, the long time behaviour of the
quantum system can also be studied from the dissipation time point of view.
The quantum setting provides a natural framework for extension of the
present work, which we will address in a separate paper \cite{FNW}.

\medskip

{\bf Acknowledgments:} This project started during a stay of S.N. at the Mathematical Sciences
Research Institute (Berkeley), the support of which is gratefully acknowledged. The authors
would like to thank Prof. S.~De~Bi\`evre and Prof. B. Nachtergaele 
for putting them together in contact.


\appendix

\section{Proofs of some elementary facts}
\label{app}


\subsection{Proof of Lemma \ref{Mom-Four}\label{app1}}

We use the following upper bound: for any $\alpha\in(0,2]$, there is a constant $C_\alpha$ 
such that
\bea
\label{cosin}
\forall x\in\IR,\quad 0\leq 1-\cos(2\pi x)\leq C_\alpha |x|^\alpha.
\eea
Besides, one has the asymptotics $1-\cos(x)\approx \frac{x^2}{2}$ for small $x$. We simply apply these
estimates to the following integral:
\bean
1-\hat g(\bxi)&=&\int_{\IR^d}(1-\cos(2\pi\bx\cdot\bxi))g(\bx)d\bx\\
&\leq&\int_{\IR^d}C_\alpha |\bx\cdot\bxi|^\alpha g(\bx)d\bx\\
&\leq& C_\alpha|\bxi|^\alpha\int_{\IR^d}|\bx|^\alpha g(\bx)d\bx=C_\alpha M_\alpha |\bxi|^\alpha.
\eean
In the case $g$ admits a second moment, we have in the limit $\bxi\to 0$:
\bean
\int_{\IR^d}(1-\cos(2\pi\bx\cdot\bxi))g(\bx)d\bx&\approx& \int_{\IR^d}2\pi^2(\bx\cdot\bxi)^2 g(\bx)d\bx\\
&\approx&2\pi^2|\bxi|^2\int_{\IR^d}(\bx\cdot\hat{\bxi})^2 g(\bx)d\bx,
\eean
where we have used the notation $\hat{\bxi}=\frac{\bxi}{|\bxi|}$ for any $\bxi\neq 0$.\hfill$\blacksquare$


\subsection{Proof of Proposition \ref{nke} \label{app2}}
The statement i) is standard in the context of distributions \cite[p.157]{Y}. 
In our case, assume that $f\in L^2$ is normalized to unity and consider an arbitrary small $\del>0$. Since
$f\in L^2(\IT^d)$, there exists $K>0$ s.t.
$\sum_{|\bk|\geq K} |\hat{f}(\bk)|^2<\del$. Since $\hat g$ is continuous and $\hat g(0)=1$, 
there exists $\eta$ such that $(1-\hat{g}(\bxi))^2<\del$ if $|\bxi|<\eta$. 
Thus using spectral decomposition \eqref{specdecompo} of $G_\ep$, we obtain for all $\ep<\frac{\eta}{K}$
\bea\label{decompo2}
\|G_{\ep}f-f\|^2= \sum_{\bk\in \IZ^d}(1-\hat{g}(\ep\bk))^2 |\hat{f}(\bk)|^{2} 
\leq \del \sum_{|\bk|<K} |\hat{f}(\bk)|^{2}+  \sum_{|\bk|>K} |\hat{f}(\bk)|^{2}
\leq 2\del.
\eea

To prove the next statement, first notice that if $g$ satisfies
the estimate \eqref{Moment-Fourier} for the exponent $\alpha$, it also satisfies it for the
exponent $\gamma\wedge\alpha$. Using once again spectral decomposition of $G_\ep$, and applying the
estimate \eqref{Moment-Fourier} with the latter exponent we get
\bea
\|G_{\ep}f-f\|^2 &\leq&\sum_{\bk\in \IZ^d}(C_{\gamma\wedge\alpha}M_{\gamma\wedge\alpha}
|\ep\bk|^{\gamma\wedge\alpha})^2 |\hat{f}(\bk)|^{2}\non\\
&\leq&(C_{\gamma\wedge\alpha}M_{\gamma\wedge\alpha})^2\ep^{2(\gamma\wedge\alpha)}
\sum_{\bk\in \IZ^d}|\bk|^{2(\gamma\wedge\alpha)}|\hat{f}(\bk)|^{2}\label{interm}\\
&\leq&(C_{\gamma\wedge\alpha}M_{\gamma\wedge\alpha})^2\ep^{2(\gamma\wedge\alpha)}
\|f\|_{H^{\gamma\wedge\alpha}}^2.\non
\eea

To obtain the last statement, we notice that any $f\in C^1(\IT^d)$ is automatically in 
$H^1(\IT^d)$, and that its gradient satisfies
\bean
\|\nabla f \|_\infty^2\geq \|\nabla f \|^2 = 4\pi^2\sum_{\bk \in \IZ^d}|\bk|^2|\hat{f}(\bk)|^2 
\geq  4\pi^2\sum_{\bk \in \IZ^d}|\bk|^{2(1\wedge\alpha)}|\hat{f}(\bk)|^2.
\eean
The inequality \eqref{interm} with $\gamma=1$ then yields the desired result.\hfill$\blacksquare$


\subsection{Proof of Corollary \ref{disto0}\label{app3}}

We prove the limit $d_\ep(1)\epto 0$ by contradiction. Assume that there is some constant $a\in(0,1)$ such that
for all $\ep>0$, the distance $d_\ep(1)>a$. We will show that the following triangle inequality holds:
\bea\label{distance}
\forall \ep>0,\quad d_\ep(1-a/2)> a/2.
\eea
First of all, notice that the assumption $d_\ep(1)>a$ means that for any $\lambda\in S^1$, $\|R_\ep(\lam)\|<a^{-1}$.
We apply the following identity \cite{Y}:
$$
R_\ep(\lam')=R_\ep(\lam)\{1+\sum_{n\geq 1}(\lam-\lam')^n R_\ep(\lam)^n\}
$$
with $\lam'=r\lam$, for $1-a<r<1$. Taking norm of both sides yields the bound
$\|R_\ep(\lam')\|\leq \frac{1}{r-(1-a)}$, uniformly w.r.t $\ep$.  
Since this upper bound holds for any $|\lam'|=r$, it shows that the spectral 
radius $r_{sp}(T_\ep)\leq 1-a$, and proves \eqref{distance} by taking $r=1-a/2$.
We can now use \eqref{distance} in the upper bound \eqref{second-upper} of 
Theorem~\ref{gb}:
this $\ep$-independent upper bound shows that ${\tau}_*$ remains finite in the 
limit $\ep\to 0$, which
contradicts Proposition~\ref{nonfinite}.\hfill$\blacksquare$


\subsection{Proof of Lemma \ref{sum-int}\label{app4}}
Considering its decay at infinity, the function $f$ is automatically in $L^2(\IR^d)$. 
The function $\hat f^2$ is the Fourier transform of the self-convolution $f\ast f$. Therefore,
using the parity of $f$ and applying the Poisson summation formula to the LHS of \eqref{sum-integ} 
yields
\bea
\label{poisson}
\ep^d\sum_{\bk\in\IZ^d}\hat f(\ep\bk)^2=\int\hat{f}^{2}(\bxi)d\bxi+ 
\sum_{0\not=\bn\in\IZ^{d}}(f\ast f)\big(\frac{\bn}{\ep}\big).
\eea
A simple computation shows that
$(f\ast f)(\bx)$ also decays as fast as $|\bx|^{-M}$. This 
piece of information is now sufficient to control the RHS of \eqref{poisson}, yielding the result,
Eq.~\eqref{sum-integ}.
\hfill$\blacksquare$

\subsection{Proof of Corollary \ref{supexp}}
It was showed in \cite{FW} that for any $\del>0$ in case of Gaussian noise
\bea
\label{norm}
\| T^{n}_{\ep}\|\leq e^{-\ep^2e^{2(1-\del)\hat{h}(F)n}}.
\eea
Using this estimate one immediately gets
\bean
C^{\ep}_{f,g}(n)=\la \bar{f}, T^{n}_{\ep}g \ra \leq \|f\|\|g\|\| T^{n}_{\ep}\| \leq \|f\|\|g\|e^{-\ep^2 \lam^n}.
\eean
Now let $f=G_{\ep}f_{0}$ and $g=G_{\ep}g_{0}$. Since the estimate (\ref{norm}) holds also in 
coarse-grained version, we have
\bean
C_{f,g}(n)&=&\la \bar{f}, U_{F}^{n}g \ra = \la G_{\ep}\bar{f_{0}}, U_{F}^{n}G_{\ep}g_{0} \ra
=\la \bar{f_{0}}, \tilde{T}^{(n)}_{\ep}g_{0} \ra \\
&\leq& \|f_{0}\|\|g_{0}\|\| \tilde{T}^{(n)}_{\ep}\| \leq \|f_{0}\|\|g_{0}\|e^{-\ep^2 \lam^n}.
\eean

\end{document}